\documentclass[graybox]{svmult}

\usepackage{type1cm}        
\usepackage{epsfig}
\usepackage{makeidx}         
\usepackage{graphicx}        
\usepackage{epstopdf}
\usepackage{multicol}        
\usepackage[bottom]{footmisc}

\usepackage{url}
\usepackage{newtxtext}       %
\usepackage{newtxmath}       
\usepackage{tikz}

\makeindex             

\def\IR{{\mathbb R}}
\def\IC{{\mathbb C}}
\def\IZ{{\mathbb Z}}

\def\IN{{\mathbb N}}

\def\IL{{\mathbb L}}

\def\IV{{\mathbb V}}
\def\IW{{\mathbb W}}
\newcommand{\sIL}{{{{\mathbb L}_s}}}

\newcommand{\bSigma}{{\bf \Sigma}}
\newcommand{\bA}{{\bf A}}
\newcommand{\bB}{{\bf B}}
\newcommand{\bC}{{\bf C}}

\newcommand{\bE}{{\bf E}}
\newcommand{\bF}{{\bf F}}
\newcommand{\bG}{{\bf G}}
\newcommand{\bS}{{\bf S}}
\newcommand{\bY}{{\bf Y}}

\newcommand{\bM}{{\bf M}}
\newcommand{\bN}{{\bf N}}
\newcommand{\bI}{{\bf I}}

\newcommand{\bP}{{\bf P}}

\newcommand{\bX}{{\bf X}}
\newcommand{\bx}{{\bf x}}
\newcommand{\by}{{\bf y}}
\newcommand{\bu}{{\bf u}}

\newcommand{\bc}{{\bf c}}

\newcommand{\bb}{{\bf b}}

\newcommand{\bz}{{\bf z}}

\newcommand{\bff}{{\bf f}}

\newcommand{\bPhi}{ \boldsymbol{\Phi} }

\newcommand{\blambda}{\boldsymbol{\lambda}}

\newcommand{\bmu}{\boldsymbol{\mu}}

\newcommand{\bLambda}{\boldsymbol{\Lambda}}

\definecolor{green}{rgb}{0,0.5,0}
\tikzstyle{format} = [draw, thin, fill=white]
\tikzstyle{format1} = [draw,white]

\begin{document}
\title*{On Bilinear Time Domain Identification {and Reduction in the Loewner Framework}}
\author{D.S. Karachalios, I.V. Gosea, A.C. Antoulas}
\institute{D.S. Karachalios \at MPI, Magdeburg, Germany, \email{karachalios@mpi-magdeburg.mpg.de}
\and I.V. Gosea \at MPI, Magdeburg, Germany, \email{gosea@mpi-magdeburg.mpg.de} \and A.C. Antoulas \at MPI, Magdeburg, Germany,\at Baylor College, Houston, USA,\at Rice University, Houston, USA, \email{aca@rice.edu}}
%
%
\maketitle

\vspace{-5mm}

\abstract*{{The \textit{Loewner framework}-(LF) in combination with \textit{Volterra series}-(VS) offers a non-intrusive approximation method that is capable of identifying bilinear models from time-domain measurements. This method uses harmonic inputs which establish a natural way for data acquisition. For the general class of nonlinear problems with VS representation, the growing exponential approach allows the derivation of the generalized kernels, namely \textit{symmetric generalized frequency response functions} - (GFRFs). In addition, the homogeneity of the Volterra operator determines the accuracy in terms of how many kernels are considered. For the weakly nonlinear setup, only a few kernels are needed to obtain a good approximation. In this direction, the proposed adaptive scheme is able to improve the estimations of the computationally non-zero kernels. The Fourier transform associates these measurements with the derived GFRFs and the LF makes the connection with system theory. In the linear case, the LF associates the so-called S-parameters with the linear transfer function by interpolating in the frequency domain. The goal of the proposed method is to extend identification to the case of bilinear systems from time-domain measurements and to approximate other general nonlinear systems (by means of the Carleman bilinearizarion scheme). By identifying the linear contribution with the LF, a considerable reduction is achieved by means of the SVD. The fitted linear system has the same McMillan degree as the original linear system. Then, the performance of the linear model is improved by augmenting a special nonlinear structure. In a nutshell, we learn reduced-dimension bilinear models directly from a potentially large-scale system that is simulated in the time domain. This is done by fitting first a linear model, and afterwards, by fitting the corresponding bilinear operator.}}

\abstract{{The \textit{Loewner framework}-(LF) in combination with \textit{Volterra series}-(VS) offers a non-intrusive approximation method that is capable of identifying bilinear models from time-domain measurements. This method uses harmonic inputs which establish a natural way for data acquisition. For the general class of nonlinear problems with VS representation, the growing exponential approach allows the derivation of the generalized kernels, namely \textit{symmetric generalized frequency response functions} - (GFRFs). In addition, the homogeneity of the Volterra operator determines the accuracy in terms of how many kernels are considered. For the weakly nonlinear setup, only a few kernels are needed to obtain a good approximation. In this direction, the proposed adaptive scheme is able to improve the estimations of the computationally non-zero kernels. The Fourier transform associates these measurements with the derived GFRFs and the LF makes the connection with system theory. In the linear case, the LF associates the so-called S-parameters with the linear transfer function by interpolating in the frequency domain. The goal of the proposed method is to extend identification to the case of bilinear systems from time-domain measurements and to approximate other general nonlinear systems (by means of the Carleman bilinearizarion scheme). By identifying the linear contribution with the LF, a considerable reduction is achieved by means of the SVD. The fitted linear system has the same McMillan degree as the original linear system. Then, the performance of the linear model is improved by augmenting a special nonlinear structure. In a nutshell, we learn reduced-dimension bilinear models directly from a potentially large-scale system that is simulated in the time domain. This is done by fitting first a linear model, and afterwards, by fitting the corresponding bilinear operator.}}

\section{Introduction}
\label{sec:1}
In natural sciences, evolutionary phenomena can be modelled as dynamical systems. An ever-increasing need for improving the approximation accuracy has motivated  including more involved and detailed features in the modelling process, thus inevitably leading to larger‐scale dynamical systems \cite{morAnt05}. To overcome this problem, efficient finite methods heavily rely on \textit{model reduction}. Model reduction methods can be classified into two broad categories, namely, \textit{SVD-based} and \textit{Krylov-based} (moment-matching).

The most prominent among the SVD-based methods is \textit{balanced truncation} (BT). In general, balancing methods are based on the computation of controllability and observability \textit{gramians} and lead to the elimination of state variables which are difficult to reach and to observe. Besides having high computational cost of solving the associated matrix Lyapunov equations, the advantages of balancing methods include the preservation of stability and an a priori computable error bound. For more details on these topics as well as on other model reduction methods not treated here (e.g., proper orthogonal decomposition (POD)/reduced basis (RB)), we refer the reader to the book \cite{morAnt05} and the surveys \cite{morBauBF14,morBenGW15}.

One way to perform model reduction is by employing, \textit{tangential interpolation}. These methods are known as rational \textit{Krylov methods} or \textit{moment-matching} methods. Krylov-based methods are numerically efficient and have lower computational cost, but in general the preservation of other properties (e.g., stability or passivity) is not automatic. For an extensive study in \textit{interpolatory model reduction}, we refer the reader to the recent book \cite{ABG20}. In what follows we will consider exclusively \textit{interpolatory model reduction methods} and, in particular the LF. For recent surveys on the LF, see \cite{morAndA90,tutorial,KGAHandbook}. The sensitivity to noise in the LF was already discussed in \cite{Lefteriu2010,drma2019learning}.

When \textit{input-output} data are offered, \textit{data-driven} methods such as the \textit{Loewner framework} - (LF), \textit{dynamic mode decomposition}-(DMD) \cite{schmid_2010}, \textit{sparse identification of nonlinear systems (with control)} - SINDYc in \cite{SINDYc}, \textit{vector fitting}-(VF) \cite{VF}, {\textit{Hankel} \cite{IsidoriMinBil}} or \textit{subspace} methods {\cite{IonitaTime,Juang,BARTEE199447}}, {\textit{moment-matching} \cite{SCARCIOTTI2017340}} and \textit{operator inference} \cite{PW2016}, remain the only feasible approaches for recovering the hidden information.

DMD-based methods represent viable alternatives that require state-derivative estimations. In such cases, the library construction allows a broad nonlinear identification. Nevertheless, the dimension of  this library could be very large when accurate discretization of the physical domain is performed.

\begin{figure}[h!]
\centering
\begin{tikzpicture}[node distance=5cm, auto, thick]
    \draw[->] node[format] (tex) {$\bu(t)$};
    \draw[->] node[format, right of=tex] (dvi) {$\bSigma:~\dot{\bx}(t)=\bff(\bx(t),\bu(t))$,~$\by(t)=\bz(\bx(t),\bu(t))$}
                  (tex) edge node {input} (dvi);
    \draw[->] node[format, right of=dvi] (ps) {$\by(t)$}
                  (dvi) edge node {output} (ps);
\end{tikzpicture}
\caption{Mathematical formalism for evolutionary phenomena.}
\label{fig: ISO}
\end{figure}
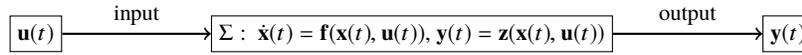

 While the underlying dynamical system acts as a black box, model identification tools are important for the reliability of the discovered models (i.e., stability, prediction). At the same time, these discovered models might have large dimension and hence are not suitable for fast numerical simulation {and control}.  The LF is a direct data-driven interpolatory method able to identify and reduce models derived directly from measurements. For measured data in the frequency domain, the LF is well established for linear and nonlinear systems (e.g., bilinear or quadratic-bilinear systems) see \cite{AGIbil,GAQuad_Bil}. In the case of time-domain data, the LF was already applied for approximating linear models \cite{IonitaTime,TimeLinearLoewner,fosong2019timedomain}. {As the aim of this paper is to extend the identification and reduction procedure to the class of bilinear systems from time-domain data, we start our analysis by introducing the mathematical description of the DAEs\footnote{DAEs: Differential Algebraic Equations.} which can simulate the \textit{input-$u(t)$} to \textit{output-$y(t)$} relation as depicted in Fig.\;\ref{fig: ISO}. The differential and algebraic operators are denoted with $\bff$ and, respectively, with $\bz$. To achieve this goal, all the important steps from nonlinear system theory and interpolatory model reduction are summarized.}
 
\subsection{Outline of the paper}
The rest of the paper is organized as follows:
\begin{itemize}
\item Section 2 contains a brief description of system theory starting from the linear case followed by extensions to the nonlinear case by means of the Volterra series representation. The single-input and single-output case is addressed for both frequency and time-domain representations.
\item Section 3 introduces the Loewner framework as an interpolatory tool for model approximation; the results that are presented here actually set the foundation for identification and reduction of linear time-invariant systems. 
\item Section 4 introduces a special class of nonlinear systems, e.g., bilinear systems. The theoretical discussion for analyzing such systems starts with the growing exponential approach and the derivation of the generalized frequency response functions (GFRFs) up to the case where a double-tone input is assumed. In addition, the kernel separation strategy for improving the measurements and the linear identification/reduction part is presented. A concise algorithm that summarizes the method is presented. 
\item Section 5 presents the numerical experiments performed in order to illustrate the practical applicability of the newly-proposed method. This section includes both a simple (low-dimensional) example as well as a large-scale example. 
\item Section 6 presents the concluding remarks and also some potential future developments of the current method.
\end{itemize}
\vspace{-5mm}
\section{System theory preliminaries} 
In this section, we will briefly present some important material from system theory starting from the linear case. 
\vspace{-5mm}
\subsection{Linear systems}\label{sec:LTI}
Consider SISO linear, time-invariant systems with $n$ internal variables (called "states" whenever the matrix $\bE$ is non-singular).
\begin{equation} \label{sysdef}
\small
\bSigma_{l}:~\left\{\begin{aligned}
\bE\,\dot\bx (t)&=\bA\bx (t)+\bb u (t),\quad\\[2mm]
y (t)&=\bc\bx (t),~t\geq 0,
\end{aligned}\right.
\end{equation}
where $\bE,~\bA\in\IR^{n\times n},~\bc\in\IR^{n\times 1},~\bc\in\IR^{1\times n}$. In the sequel we will restrict our attention to invertible matrix $\bE$ and with zero d-term ($d=0$) in the state-output equation\footnote{The state-output  equation often is represented as $y(t)=\bc\bx(t)+du(t)$.}. 
The explicit solution with the \textit{convolution integral}\footnote{$(h\ast u)(t)=\int_{-\infty}^{\infty}h(\tau)u(t-\tau)d\tau$.} notation and the time-domain linear kernel $h(t)$ as the \textit{impulse response} of the system can be written as:
\begin{equation}
\small
y(t)=\bc e^{\bA t}\bx (0)+(h\ast u)(t),~t\geq 0.
\end{equation} 
By assuming zero initial conditions and performing a \textit{Laplace transform}, we obtain the transfer function description:
\begin{equation}
\small
H(s)=\frac{Y(s)}{U(s)}=\bc(s\bI-\bA)^{-1}\bb,~s\in\IC,
\end{equation}
where $Y(s),U(s)$ stand for the \textit{input} and the \textit{output} in the frequency domain.  
\subsection{Nonlinear systems}
A large class of nonlinear systems can be described by means of the
Volterra-Wiener approach in \cite{Rug81}. Other relevant works on nonlinear systems and nonlinear modelling/identification include Schetzen (1980), Chen and Billings (1989), Boyd and Chua (1985) et. al. 

\begin{figure}[h!]
\centering
\begin{tikzpicture}[node distance=4cm, auto, thick]
    \draw[->] node[format] (tex) {$u(t)$};
    \draw[->] node[format, right of=tex] (dvi) {$\bSigma:$ unknown}
                  (tex) edge node {input} (dvi);
    \draw[->] node[format, right of=dvi] (ps) {$\by(t)$}
                  (dvi) edge node {output} (ps);
\end{tikzpicture}
\caption{The input-output mapping from the data-driven perspective with the unknown system $\bSigma$. Specific structures of the unknown system can be assumed/inspired by the physical problem. For instance, if the underlying physical phenomenon is fluid flow inside a control volume, quadratic models should be constructed e.g., \cite{GAQuad_Bil}.}
\label{fig: ISO2}
\end{figure}
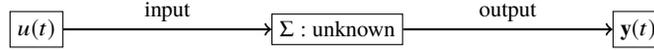

The aim in this study is to identify and reduce special types of nonlinear systems (s.a., bilinear) from time-domain measurements. By knowing only the input and the simulated or measured output in the time domain as in Fig. \ref{fig: ISO2}, we will identify the hidden model. In such situations where only snapshots are available, beyond the linear fit which is well established a nonlinear fit of a special type will be developed. 
\subsubsection{Approximation of nonlinear systems (Volterra series)}
The \textit{input-output} relationship for a wide class of nonlinear systems \cite{Rug81} can be approximated by a Volterra series for sufficiently high $N$ as:
\begin{equation}\label{eq:VoltOutput}
\small
y(t)=\sum_{n=1}^{N}y_{n}(t),~~y_{n}(t)=\int_{-\infty}^{\infty}\cdots\int_{-\infty}^{\infty}h_{n}(\tau_{1},\cdots,\tau_{n})\prod_{i=1}^{n}u(t-\tau_{i})d\tau_{i},
\end{equation}
where $h_{n}(\tau_{1},\cdots,\tau_{n})$ is a \textit{real-valued} function of $\tau_{1},
\cdots,\tau_{n}$ known as the $n^{th}$ order Volterra kernel.
\begin{definition}\label{eq:nthVolterraKernel}
The $n^{th}$ order generalized frequency response function (GFRF) is defined as:
\begin{equation}\label{eq:nthVolterraKernelH}
\small
H_{n}(j\omega_{1},\cdots,
j\omega_{n})=\int_{-\infty}^{\infty}\cdots\int_{-\infty}^{\infty}h_{n}(\tau_{1},\cdots,\tau_{n})e^{\left(-j\sum_{i=1}^{n}\omega_{i}\tau_{i}\right)}d\tau_{1}\cdots d\tau_{n},~(j^2=-1)
\end{equation}
which is the multidimensional Fourier\footnote{As the frequency $s=j\omega$ lies on the imaginary axis, the Laplace transform simplifies in most cases to Fourier transform (e.g. for square integrable functions).} transform of $h_{n}(\tau_{1},\cdots,\tau_{n})$.
\end{definition}
By applying the inverse Fourier transform of the $n^{th}$ order GFRF, Eq.\;\eqref{eq:nthVolterraKernelH} can be written as:
\begin{equation}
\small
y_{n}(t)=\frac{1}{{(2\pi)}^n}\int_{-\infty}^{\infty}\cdots\int_{-\infty}^{\infty}H_{n}(j\omega_{1},\cdots,
j\omega_{n})\prod_{i=1}^{n}U(j\omega_{i})e^{j(\omega_{1}+\cdots+\omega_{n})t}d\omega_{1}\cdots\omega_{n}.
\end{equation}
The $n^{th}$ Volterra operator is defined as:
\begin{equation}
\small
V_{n}(u_1,u_2,...,u_n)=\int_{-\infty}^{\infty}\cdots\int_{-\infty}^{\infty}h_{n}(\tau_{1},...,\tau_{n})\prod_{i=1}^{n}u_{i}(t-\tau_{i})d\tau_{i},
\end{equation}
so that $y_{n}=V_{n}(u,u,...,u)$ holds true.
\begin{important}{Homogeneity of the Volterra operator} The map $u(t)\rightarrow y_{n}(t)$ is homogeneous of degree $n$, that is, $\alpha u\rightarrow \alpha^n y_{n},~\alpha\in\IC$. Each Volterra kernel $h_{n}(t)$ determines a symmetric multi-linear operator. Small amplitudes (e.g., $|\alpha|<\epsilon$) will allow ordering the nonlinear terms in such a way that terms with large powers of the amplitude ($\alpha^n$) will be negligible. That is precisely the sense of approximating weakly nonlinear systems with Volterra series.
\end{important}
\subsubsection{A single-tone input}
Consider the excitation of a system with an input consisting of two complex exponentials as in Eq.\;\eqref{eq:Acos}. Such inputs are typically used in chemical engineering applications as \cite{PDM}.
\begin{equation}\label{eq:Acos}
\small
u(t)=A\cos(\omega t)=\left(\frac{A}{2}\right)e^{j\omega t}+\left(\frac{A}{2}\right)e^{-j\omega t}.
\end{equation}
By using the above input in Eq.\;\eqref{eq:VoltOutput}, we can derive the first Volterra term with $n=1$ as:
\begin{equation}
\small
\begin{aligned}
y_{1}(t)&=\int_{-\infty}^{\infty}h_{1}(\tau_{1})[u(t-\tau_{1})]d\tau_{1}\\
        &=\frac{A}{2}e^{j\omega t}\underbrace{\int_{-\infty}^{\infty}h_{1}(\tau_{1})e^{-j\omega \tau_{1}}d\tau_{1}}_{H_{1}(j\omega)}+\frac{A}{2}e^{-j\omega t}\underbrace{\int_{-\infty}^{\infty}h_{1}(\tau_{1})e^{j\omega\tau_{1}}d\tau_{1}}_{H_{1}(-j\omega)}\Rightarrow\\
        y_{1}(t)&=\frac{A}{2}\bigg(e^{j\omega t}H_{1}(j\omega)+e^{-j\omega t}H_{1}(-j\omega)\bigg).
\end{aligned}
\end{equation}
Similarly, for the 2nd term we can derive:
\begin{equation}
\small
y_{2}(t)=\bigg(\frac{A}{2}\bigg)^2\bigg[e^{2j\omega}H_{2}(j\omega,j\omega)+2e^{0}H_{2}(j\omega,-j\omega)+e^{-2j\omega}H_{2}(-j\omega,-j\omega)\bigg].
\end{equation}
\begin{remark}{(Conjugate symmetry)}: $H_{2}^*(j\omega,-j\omega)=H_{2}(-j\omega,j\omega),~\forall\omega\in\IR$.
 \end{remark}
The input amplitude is $A$, the angular frequency is $\omega$, the imaginary unit is $\mathrm{j}$, the first order response function is $H_{1}(j\omega)$, and $H_{n}(j\omega,...,j\omega)$, for $n\geq 2$, are the higher-order FRFs or GFRFs. Then, the $n^{th}$ Volterra term can be written as:
\begin{equation}\label{eq:nthHarm}
\small
y_{n}(t)=\left(\frac{A}{2}\right)^{n}\sum_{p+q=n}{}^{n}C_{q}H_{n}^{p,q}(j\omega)e^{j\omega_{p,q}t},~\omega_{p,q}=(p-q)\omega.
\end{equation}
where the following notations have been used: 
\begin{equation}\label{eq:notation}
\small
H_{n}^{p,q}(j\omega)=H_{n}(\underbrace{j\omega,...,j\omega}_{p-times};\underbrace{-j\omega,...,-j\omega}_{q-times}),~\omega_{p,q}=(p-q)\omega,~{}^{n}C_{q}=\frac{n!}{q!(n-q)!}.
\end{equation}
\subsubsection{Time domain representation of harmonics}
The $m^{th}$ harmonic in the \textit{time domain} can be computed by collecting the identical exponential power coefficients from Eq.\;\eqref{eq:ntimeharm} and by setting $p-q=m$, with $p=m+i-1$ and $q=i-1$ in Eq.\;\eqref{eq:nthHarm}. Hence, it follows that:
\begin{equation}\label{eq:ntimeharm}
\small
y_{m^{th}}(t)=\sum_{i=1}^{\infty}\left(\frac{A}{2}\right)^{m+2i-2}{}^{m+2i-2}C_{i-1}H_{m+2i-2}^{m+i-1,i-1}(j\omega)e^{jm\omega t}.
\end{equation}
\subsubsection{Frequency domain representation of harmonics}
\begin{figure}[h!]
\sidecaption[t]
\centering
\includegraphics[width=75mm,height=33mm]{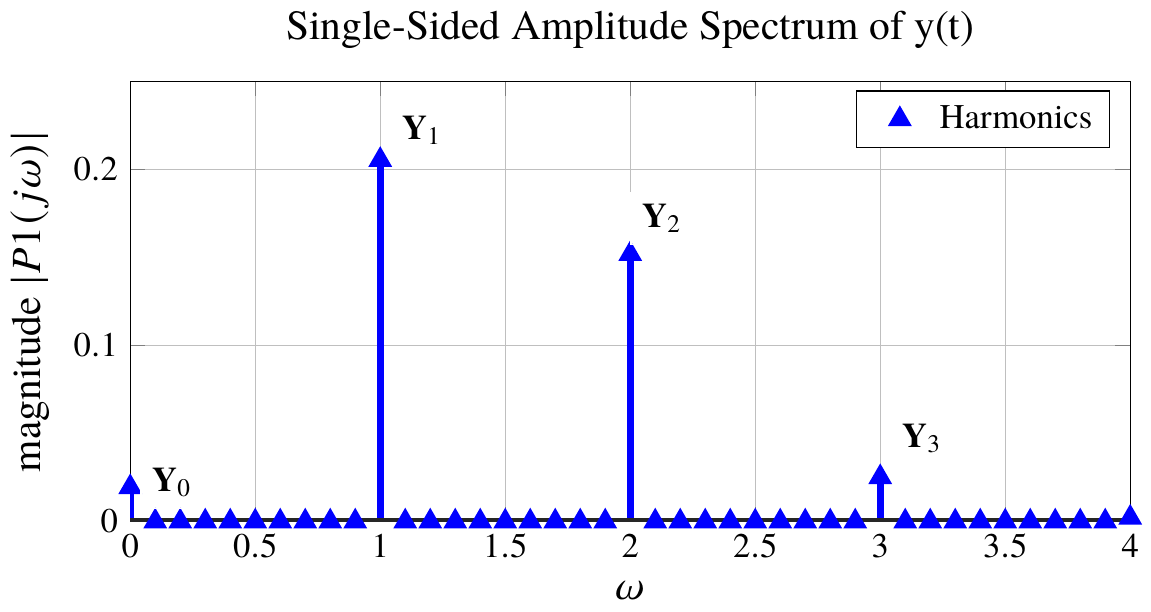}
\caption{An instance of the single-sided power spectrum with a singleton input with $\omega=1$ is depicted. The underlying system is nonlinear and as a result higher harmonics appeared with a DC (Direct Current - non-periodic) term as well.}
\label{fig:Spectrum}
\end{figure}

The $m^{th}$ harmonic in the \textit{frequency domain} by applying single-sided Fourier transform in Eq.\;\eqref{eq:ntimeharm} is the following: 
\begin{equation}\label{eq:fourierOutput}
\small
Y_{m^{th}}(jm\omega)=\sum_{i=1}^{\infty}2\left(\frac{A}{2}\right)^{m+2i-2}{}^{m+2i-2}C_{i-1}H_{m+2i-2}^{m+i-1,i-1}(j\omega)\delta(jm\omega).
\end{equation}
where $\delta(\cdot)$ is the Dirac delta distribution. When a single-tone input excites a nonlinear dynamical system, the steady state frequency response is characterized by a spectrum with higher harmonics (as can be seen, for example, in Fig.\;\ref{fig:Spectrum}). This behavior is not observed in the linear case, where only one harmonic appears at the input frequency.
\vspace{-5mm}
\section{The Loewner framework}
We start with an account of the Loewner framework (LF) in the linear case \cite{morAndA90,tutorial,KGAHandbook}. The (LF) is an interpolatory method that seeks reduced models whose transfer function matches that of the original system at selected interpolation points. An important attribute is that it provides a trade-off between accuracy of fit and complexity of the model. It constructs models from given frequency data in a straightforward manner. In the case of SISO systems, we have the rational scalar interpolation problem to solve. 

Consider a given set of {complex data} as:
$$
\{\left(s_k,f_k(s_k)\right)\in\IC\times\IC:k=1,\ldots,2n)\}.
$$
We partition the data in two disjoint sets:
$$
\bS=[\underbrace{s_1,\ldots,s_n}_{\bmu},\underbrace{s_{n+1},\ldots,s_{2n}}_{\blambda}],~\bF=[\underbrace{f_1,\ldots,f_n}_{\IV},\underbrace{f_{n+1},\ldots,f_{2n}}_{\IW}],
$$
where $\mu_i=s_i$, $\lambda_i=s_{n+i}$, $v_{i}=f_i$, $w_{i}=f_{n+i}$ for $i=1,\ldots,n$.\\[1mm]
The objective is to find $H(s)\in\IC$ such that:
\begin{equation}\label{prob:inter}
\small
H(\mu_i)=v_{i},~i=1,\ldots,n,~\text{and}~H(\lambda_j)=w_{j},~j=1,\ldots,n.
\end{equation}
The {\it left data set} is denoted as:
\begin{equation}\label{leftdata}
\small
\bM=\left[\mu_1,\cdots,\mu_n\right]\in\IC^{1\times n},~~\IV=\left[v_1,\cdots,v_n\right]^{T}\in\IC^{n\times 1} ,
\end{equation}
while the {\it right data set} as:
\begin{equation}\label{rightdata}
\small
\bLambda=\left[\lambda_1,\cdots,\lambda_n\right]^{T}\in\IC^{n\times 1},~~\IW=[w_1,\cdots, w_n]\in\IC^{1\times n}.
\end{equation}
Interpolation points are determined by the problem or are selected to achieve given model reduction goals. For ways of choosing the interpolation grids and of partitioning the data into the left and right sets, we refer the reader to the recent survey \cite{KGAHandbook}.  
\subsection{The Loewner matrix}
Given a row array of complex numbers
$(\mu_j,v_j)$, $j=1,\ldots,{n}$, and
a column array, $(\lambda_i,w_i)$, $i=1,\ldots,{n},$  
(with $\lambda_i$ and the $\mu_j$ mutually distinct)
the associated {\it Loewner matrix} $\IL$ and the shifted {\it Loewner matrix} $\sIL$ are defined as:
$$
\hspace*{-4mm}
\IL\!=\!\left[\!\begin{array}{ccc}
\frac{v_1-w_1}{\mu_1-\lambda_1} & 
\cdots &
\frac{v_1-w_{n}}{\mu_1-\lambda_{n}} \\
\vdots & \ddots & \vdots \\
\frac{v_n-w_1}{\mu_n-\lambda_1} & 
\cdots &
\frac{v_n-w_{n}}{\mu_n-\lambda_{n}} 
\end{array}\!\right]\!\in\!\IC^{n\times n},~
\sIL\!=\!\left[\!\begin{array}{ccc}
\frac{\mu_1v_1-\lambda_1w_1}{\mu_1-\lambda_1} & 
\cdots &
\frac{\mu_1v_1-\lambda_{n}w_{n}}{\mu_1-\lambda_{n}} \\
\vdots & \ddots & \vdots \\
\frac{\mu_n v_n-\lambda_1 w_1}{\mu_n-\lambda_1} & 
\cdots &
\frac{\mu_n v_n-\lambda_{n}w_{n}}{\mu_n-\lambda_{n}} 
\end{array}\!\right]\!\in\!\IC^{n\times n}.
$$  
\begin{definition}\label{def:McMillan}
If $g$ is rational, i.e., $g(s)=\frac{p(s)}{q(s)}$, for appropriate polynomials $p$, $q$, the McMillan degree or the complexity of $g$ is~$\mbox{deg}\,g=\max\{\mbox{deg}(p),\mbox{deg}(q)\}$.
\end{definition}

Now, if $w_i=g(\lambda_i)$, and $v_j=g(\mu_j)$, are 
\textit{samples} of a rational function $g$, the \textit{main property} 
of Loewner matrices asserts the following.
\begin{theorem}\cite{morAndA90}
Let $\IL$ be as above. If $
k,q\geq {\rm deg}\,g$, then $\,{\rm rank}\, \IL = {\deg}\,g$.

In other words the rank of $\IL$ encodes the complexity of the underlying 
rational function $g$. Furthermore, the same result holds for
matrix-valued functions $g$.
\end{theorem}
\subsection{Construction of interpolants}
If the pencil ~$(\sIL,\,\IL)$~ is regular, then ~
$\bE=-\IL,~~ \bA=-\sIL,~~ \bb=\IV,~~ \bc=\IW$,~ 
is a minimal realization of an interpolant for the data, i.e., 
$H(s)=\IW(\sIL-s\IL)^{-1}\IV$. Otherwise, as shown in \cite{morAndA90}, the problem in Eq.\;\eqref{prob:inter} has a solution provided that 
\begin{equation*}\label{assumption}
\small
\hspace*{-2mm}
\mbox{rank}\,\left[s\,\IL-\sIL\right]=\mbox{rank}\,\left[\IL,\ \, \sIL\right]=
\mbox{rank}\,\left[\!\begin{array}{c}\IL\\\sIL\end{array}\!\right]\!= {r},
\end{equation*}
for all ~$s\in\{\mu_i\}\cup\{\lambda_j\}$.~ 
Consider then the thin SVDs: 
\begin{equation*}\label{prmat}
\small
\left[\IL,\ \, \sIL\right]=\bY\widehat{\Sigma}_{ {r}}\tilde{\bX}^*,~~
\left[\begin{array}{c}\IL\\\sIL\end{array}\right] = {\tilde\bY}\Sigma_{ {r}} \bX^*,
\end{equation*}
where ~$\widehat{\Sigma}_{ {r}}$, $\Sigma_{ {r}}$ $\in$ $\IR^{{ {r}}\times  {r}}$,~
$\bY \in\IC^{n\times  {r}}$,
$\bX$ $\in$ $\IC^{n \times  {r}}$, $\tilde{\bY} \in\IC^{2n\times  {r}}$, $\tilde{\bX}$ $\in$ $\IC^{r \times  {2n}}$.  

\begin{remark}
{$r$} can be chosen as the {\it numerical rank} 
(as opposed to the {\it exact rank}) of the Loewner pencil. 
\end{remark} 

\begin{theorem}
The quadruple $(\tilde{\bA},\tilde{\bb},\tilde{\bc},\tilde{\bE})$ of size 
~$ {r}\times  {r}$, ~$ {r}\times  {r}$, ~$ {r}\times 1$, ~$1\times  {r}$,~ given by:
\begin{equation*} \label{redundant}
\small
\tilde{\bE}  = -\bY^T\IL \bX ,~~\tilde{\bA}  = -\bY^T\sIL \bX, ~~\tilde{\bb}  = \bY^T\IV,~~
\tilde{\bc}  = \IW \bX ,
\end{equation*}
is a descriptor realization of an (approximate) interpolant of the data with McMillan degree $r=rank(\IL)$, where $\tilde{H}(s)=\tilde{\bc}(s\tilde{\bE}-\tilde{\bA})^{-1}\tilde{\bb}$.
\end{theorem}

For more details on the construction/identification of linear systems with the (LF), we refer the reader to \cite{tutorial,KGAHandbook} where both the SISO and MIMO cases are addressed together with other more technical aspects (e.g., how to impose the construction of real-valued models, etc.).
\vspace{-5mm}
\section{The special case of bilinear systems}
In recent years, projection-based Krylov methods have extensively been applied for model reduction of bilinear systems. We mention the following contributions \cite{AGIbil,morAhmBB17,morBaiS06,morBenB12b,morBenBD11,morBenD11,morBre13,morFlaG15,morPhi03} and the references within. 

Scalar bilinear systems are described by the set of matrices; $\bSigma_{b}=(\bA,\bN,\bb,\bc,\bE)$ and characterized by the following equations:
\begin{equation}\label{eq:bilinear}
\bSigma_{b}: \left\{\begin{aligned}
\bE\dot{\bx}(t)&=\bA\bx(t)+\bN\bx(t)u(t)+\bb u(t),\\
          y(t)&=\bc\bx(t),
\end{aligned}\right.
\end{equation}
where $\bE,\bA,\bN\in\IR^{n\times n}$, $\bb\in\IR^{n\times 1}$, $\bc\in\IR^{1\times n}$, and $\bx\in\IR^{n\times 1},u,y\in\IR$. In what follows, we restrict our analysis to systems with non-singular $\bE$ matrices (e.g. identity matrix).
\subsection{The growing exponential approach}
The properties of the growing exponential approach can be adapted readily to the problem of \textit{finding transfer functions} for constant-parameter (stationary) state equations. Let us consider the bilinear model in Eq.\;\eqref{eq:bilinear} with zero initial conditions. A single-tone input with amplitude $A<1$ is considered as in Eq.\;\eqref{eq:Acos}. 
\begin{equation}
\small
u(t)=A\cos(\omega t)=\frac{A}{2}e^{j\omega t}+\frac{A}{2}e^{-j\omega t}=a e^{j\omega t}+a e^{-j\omega t},
\end{equation}
where $a=A/2$ and $a\in(0,\epsilon)$ with $0<\epsilon<1/2$ and for all $t\geq  0$. The steady state solution for the differential equation in Eq.\;\eqref{eq:bilinear} can be written as:
\begin{equation}
\small
\bx(t)=\sum_{p,q\in\IN}^{\infty}\bG_{n}^{p,q}(\underbrace{j\omega,\ldots,j\omega}_{p-times},\underbrace{-j\omega,...,-j\omega}_{q-times})a^{p+q}e^{j\omega(p-q)t}
\end{equation}
The symbol\footnote{$\bG_{n}^{p,q}=\bG(\underbrace{j\omega,...,j\omega}_{p-times};\underbrace{-j\omega,...,-j\omega}_{q-times})$.} $\bG_{n}^{p,q}$ denotes the $n^{th}$ input to state frequency response containing $p$-times the frequency $\omega$ and $q$-times the frequency $-\omega$. By substituting in Eq.\;\eqref{eq:bilinear} and collecting the terms of the same exponential (as the $e^{j\omega_m t}$), we can derive the input to state frequency responses $\bG_n$ for every $n$ as follows:
\begin{equation*}
\small
\begin{aligned}
&\bE\dot{\bx}(t)-\bA\bx(t)=\bN\bx(t)u(t)+\bb u(t)\Rightarrow\sum_{p,q\in\IN}^{\infty}\left(j\omega(p-q)\bE-\bA\right)\bG_{n}^{p,q}a^{p+q}e^{j\omega(p-q)t}=\\
=&\bN\left(\sum_{p,q\in\IN}^{\infty}\bG_{n}^{p,q}a^{p+q+1}e^{j\omega(p+1-q)t}+\sum_{p,q\in\IN}^{\infty}\bG_{n}^{p,q}a^{p+q+1}e^{j\omega(p-q-1)t}\right)+\bb(ae^{j\omega t}+ae^{-j\omega t})
\end{aligned}
\end{equation*}
For the first choices of $p$ and $q$ up to $p+q\leq 2,~(1,0),(0,1),(2,0),(0,2),(1,1)$ and by denoting the resolvent $\bPhi(j\omega)=\left(j\omega\bE-\bA\right)^{-1}\in\IC^{n\times n}$, c.t. conjugate terms, we derive the first set of terms:
\begin{equation*}
\small
\begin{aligned}
&\bPhi(j\omega)^{-1}\bG_{1}^{1,0}ae^{j\omega t}+\bPhi(2j\omega)^{-1}\bG_{2}^{2,0}a^2e^{2j\omega t}+\bPhi(0)^{-1}\bG_{2}^{1,1}a^2+c.t.+\cdots=\\
&\bN\bG_{1}^{1,0}a^2e^{2j\omega t}+\bN\bG_{2}^{2,0}a^3e^{3j\omega t}+\bN\bG_{2}^{1,1}a^3e^{j\omega t}+c.t.+\cdots+\bb ae^{j\omega t}+c.t.
\end{aligned}
\end{equation*}
Collecting the same powers in both exponential and polynomial magnitudes, we compute the first and the second time/input invariant GFRFs:
\begin{equation}
\small
\begin{aligned}
\bG_{1}^{1,0}(j\omega)&=\bPhi(j\omega)\bb,\\
\bG_{2}^{2,0}(j\omega)&=\bPhi(2j\omega)\bN\bG_{1}^{1,0}=\bPhi(2j\omega)\bN\bPhi(j\omega)\bb.
\end{aligned}
\end{equation} 
Then, the following input to state transfer functions $\bG_n$ using induction are:
\begin{equation}
\small
\begin{aligned}
&\bG_{n}^{n,0}(j\omega)=\bPhi(nj\omega)\bN\bPhi((n-1)j\omega)\bN\cdots\bN\bPhi(j\omega)\bb,\\
&\bG_{n}^{0,n}(j\omega)=\bPhi(-nj\omega)\bN\bPhi(-(n-1)j\omega)\bN\cdots\bN\bPhi(-j\omega)\bb,\\
&\bG_{n}^{p,q}(j\omega)=\bPhi((p-q)j\omega)\bN\left[\bG_{n-1}^{p,q-1}(j\omega)+\bG_{n-1}^{p-1,q}(j\omega)\right],~p,q\geq 1,
\end{aligned}
\end{equation}
for $n\geq 1$ and $p+q=n$. By multiplying with the output vector $\bc$, we can further derive the input-output generalized frequency responses GFRFs as:
\begin{equation}
\small
\begin{aligned}
&H_{n}^{n,0}(j\omega)=\bc\bPhi(nj\omega)\bN\bPhi((n-1)j\omega)\bN\cdots\bN\bPhi(j\omega)\bb,\\
&H_{n}^{0,n}(j\omega)=\bc\bPhi(-nj\omega)\bN\bPhi(-(n-1)j\omega)\bN\cdots\bN\bPhi(-j\omega)\bb,\\
&H_{n}^{p,q}(j\omega)=\bc\bPhi((p-q)j\omega)\bN\left[\bG_{n-1}^{p,q-1}(j\omega)+\bG_{n-1}^{p-1,q}(j\omega)\right],~p,q\geq 1.
\end{aligned}
\end{equation}

At this point, we can write the Volterra series by using the above specific structure of the GFRFs that were derived with the growing exponential approach for the bilinear case. An important property to notice is that the $n^{th}$ kernel is a multivariate function of order $n$. It is obvious that the identification of the $n^{th}$ order FRF involves an $n$-dimensional frequency space. For that reason, next, we derive the general 2nd symmetric kernel for the bilinear case with a double-tone input. Consider:
\begin{equation}\label{eq:doubleinput}
\small
u(t)=A_{1}\cos(\omega_1 t)+A_{2}\cos(\omega_2 t)=\sum_{i=1}^{2}\alpha_{i}(e^{j\omega_i t}+e^{-j\omega_i t}),
\end{equation} 
where $\alpha_1=\frac{A_1}{2}$ and $\alpha_2=\frac{A_2}{2}$. In that case, with the growing exponential approach the state solution in steady state is
\begin{equation}
\small
\bx(t)=\sum_{m_{1},\ldots,m_{4}\in\IN}^{\infty}\bG_{n}^{m_1,m_2,m_3,m_4}\alpha_1^{m_1+m_2}\alpha_2^{m_3+m_4}e^{j((m_1-m_2)\omega_1+(m_3-m_4)\omega_2)t}.
\end{equation}

We are looking for the input to state frequency response $\bG(j\omega_1,j\omega_2)$. By substituting to the bilinear model in Eq.\;\eqref{eq:bilinear} and collecting the appropriate terms while at the same time using the symmetry $\bG(j\omega_1,j\omega_2)=\bG(j\omega_2,j\omega_1)$, we conclude that:
\begin{equation}
\small
\bG_{2}(j\omega_1,j\omega_2)=\frac{1}{2}\left[(j\omega_1+j\omega_2)\bE-\bA\right]^{-1}\bN\left[\left(j\omega_1\bE-\bA\right)^{-1}\bb+\left(j\omega_2\bE-\bA\right)^{-1}\bb\right],
\end{equation}
where by using the resolvent notation and multiplying with $\bc$, we derive the \textit{2nd order symmetric generalized frequency response function as:}
\begin{equation}
\small
H_{2}(j\omega_1,j\omega_2)=\frac{1}{2}\bc\bPhi(j\omega_1+j\omega_2)\bN\left[\bPhi(j\omega_1)\bb+\bPhi(j\omega_2)\bb\right].
\end{equation} 
\subsection{The kernel separation method}\label{sec: separation}
One way to deduce Volterra kernels is by means of interpolation. This problem is equivalent to that of estimating a polynomial with noisy coefficients. This interpolation scheme builds a linear system with a Vandermonde matrix which is invertible since the amplitudes are distinct and non-zero. The inverse of a Vandermonde matrix can be explicitly computed and there are stable ways to solve these equations \cite{Boyd1983MeasuringVK}. The recently proposed method presented in \cite{brubeck2019vandermonde} solves the exponentially ill-condition problem of the Vandermonde matrix with Arnoldi orthogonalization. The $m^{th}$ harmonic in the frequency domain is derived by applying a (single-sided)  Fourier transform. More precisely, the explicit formulation is as follows: 
\begin{equation}\label{eq:fourierOutput}
\small
\begin{aligned}
Y_{m^{th}}(jm\omega)&=\sum_{i=1}^{\infty}\underbrace{2\left(\frac{A}{2}\right)^{m+2i-2}{}^{m+2i-2}C_{i-1}}_{\alpha^{m+2i-2}}H_{m+2i-2}^{m+i-1,i-1}(j\omega)\delta(jm\omega)\\
&=\sum_{i=1}^{\infty}\alpha^{m+2i-2}H_{m+2(i-1)}^{m+i-1,i-1}(j\omega)\delta(jm\omega).
\end{aligned}
\end{equation}

We simplify the notation in order to reveal the adaptive method that will help us to estimate the GFRFs up to a specific order. Next, write the linear system of equations that connects the harmonic information with the higher Volterra kernels as follows:
\begin{equation}\label{eq:system_YMP}
\small
\begin{aligned}
\underbrace{\left[\begin{array}{c} Y_{0}(0j\omega) \\[1mm] Y_{1}(1j\omega) \\[1mm]Y_{2}(2j\omega) \\[1mm]Y_{3}(3j\omega) \\[1mm] \vdots \\[1mm] Y_{m}(mj\omega)\end{array}\right]}_{\bY_{(\alpha,\omega)}}=\underbrace{\left[\begin{array}{cccc} {\alpha^{0}} & {\alpha^{2}} & \alpha^{4} & \dots \\[1mm] {\alpha^{1}}& {\alpha^{3}} & \alpha^{5} & \dots  \\[1mm] {\alpha^{2}}& \alpha^{4} & \alpha^{6} & \dots \\[1mm]{\alpha^{3}}& \alpha^{5} & \alpha^{7} & \dots \\[1mm] \vdots & \vdots & \vdots & \vdots \\[1mm] \alpha^{m}& \alpha^{m+2} & \alpha^{m+4} & \dots  \end{array}\right]}_{\bM_{\alpha}} &\bigcirc \underbrace{\left[\begin{array}{cccc} H_{0}^{0,0}& H_{2}^{1,1} & H_{4}^{2,2} & \dots \\[1mm] H_{1}^{1,0}& H_{3}^{2,1} & H_{5}^{3,2} &\dots  \\[1mm] H_{2}^{2,0}& H_{4}^{3,1} & H_{6}^{4,2} & \dots  \\[1mm]H_{3}^{3,0}& H_{5}^{4,1} & H_{7}^{5,2} & \dots  \\[1mm] \vdots & \vdots & \vdots & \vdots \\ H_{n}^{n,0}& H_{n+2}^{n+1,1} & H_{n+4}^{n+2,2} & \dots \end{array}\right]}_{\bP_{\alpha}}\underbrace{\left[\begin{array}{c} 1 \\[1mm] 1\\[1mm] 1 \\[1mm]1 \\[1mm] \vdots \\[1mm] 1  \end{array}\right]}_{\textbf{e}_{n+1,1}}.
\end{aligned}
\end{equation}

 By introducing the Hadamard product notation\footnote{the Hadamard product is denoted with "$\circ$"; the matrix multiplication is performed element-wise.} and by substituting the $\delta$'s with ones, we can compactly rewrite the above system in the following form:
\begin{equation}\label{eq:amplitudeImprovement}
\bY_{(\alpha,\omega)}=\left[\bM_{\alpha}\circ\bP_{\omega}\right]\cdot\textbf{e}_{n+1,1}.
\end{equation}

The above system offers the level of approximation we want to achieve. Note that the frequency response $\bY$ depends on both the amplitude and the frequency, while the right hand side of Eq.\;\eqref{eq:amplitudeImprovement} reveals the separation of the aforementioned quantities. As we neglect higher order Volterra kernels, the measurement set tends to be corrupted by noise. 

\begin{important}{Kernel separation and stage $\ell$- approximation} For a given system, the procedure consists in exciting it with a single-tone input. By varying the driving frequency, as well as the amplitude, we can approximate the GFRFs by minimizing the (2-norm) of the remaining over-determined systems.
\begin{equation}
\small
\bY_{m+1,\ell}(jm\omega,\alpha_{\ell})=\left[\bM_{m+1,\ell}(\alpha_{\ell})\circ\bP_{m+1,\ell}(jm\omega)\right]\cdot\textbf{e}_{n+1,1}. 
\end{equation}
The $m$-"direction" gives us the threshold up to the specific harmonic that we measure while the $\ell$-"direction" gives us the level of the kernel separation that we want to achieve. For instance, for the 2nd stage approximation, it holds: $\ell=2$ with $\alpha^m\approx 0,~\forall m>\ell=2$.
\end{important}
\subsection{Identification of the matrix $\bN$}
The difference between linear and bilinear models is the presence of the product between the input and the state that is scaled by the matrix $\bN$. As the (LF) is able to identify the linear part ($\bA,\bb,\bc,\bE$) of the bilinear model the only thing that remains is the identification of the matrix $\bN$. The matrix $\bN$ enters linearly in the following kernels:
\begin{itemize}
\item With a single-tone input the $H_{2}^{1,1}$ can be written as: 
\begin{equation}\label{eq:2ndkernelPerDiag}
\small
H_{2}(j\omega_1,-j\omega_1)=\frac{1}{2}\bc\left(-\bA\right)^{-1}\bN\left((j\omega_1\bI-\bA)^{-1}\bb+(-j\omega_1\bI-\bA)^{-1}\bb\right).
\end{equation}
and the $H_{2}^{2,0}$ as:
\begin{equation}\label{eq:2ndkernelDiag}
\small
H_{2}(j\omega_1,j\omega_1)=\bc\left(2j\omega_1\bI-\bA\right)^{-1}\bN(j\omega_1\bI-\bA)^{-1}\bb.
\end{equation}
\item While with a double-tone input the general $H_{2}$ can be written as:
 \begin{equation}\label{eq:2ndkernel}
\small
H_{2}(j\omega_1,j\omega_2)=\frac{1}{2}\bc\left(j(\omega_1+j\omega_2)\bI-\bA\right)^{-1}\bN\left((j\omega_1\bI-\bA)^{-1}\bb+(j\omega_2\bI-\bA)^{-1}\bb\right).
\end{equation}
\end{itemize}
We introduce the following notation:
\begin{equation}
\small
\begin{aligned}
\mathcal{O}(j\omega_1,j\omega_2)&=\frac{1}{2}\bc\left(j(\omega_1+j\omega_2)\bI-\bA\right)^{-1}\in\IC^{1\times n},\\
\mathcal{R}(j\omega_1,j\omega_2)&=\left((j\omega_1\bI-\bA)^{-1}\bb+(j\omega_2\bI-\bA)^{-1}\bb\right)\in\IC^{n\times 1}.
\end{aligned}
\end{equation}
Then, Eq.\;\eqref{eq:2ndkernel} can be compactly rewritten as:
\begin{equation}
\small
H_{2}(j\omega_1,j\omega_2)=\mathcal{O}(j\omega_1,j\omega_2)\bN\mathcal{R}(j\omega_1,j\omega_2).
\end{equation}
Assume that $k$ measurements of the function $H_2$ are available (measured) for $k$ different pairs $(\omega_1,\omega_2)$. By vectorizing in respect to the measurement set, we have:
\begin{equation*}
\small
\text{The $k^{th}$ measurement}\rightarrow\underbrace{H_{2}(j\omega_{1,k},j\omega_{2,k})}_{\bY^{(k)}}=\underbrace{\mathcal{O}(j\omega_{1,k},j\omega_{2,k})}_{\mathcal{O}_{1,n}^{(k)}}\underbrace{\bN}_{n\times n}\underbrace{\mathcal{R}(j\omega_{1,k},j\omega_{2,k})}_{\mathcal{R}_{n,1}^{(k)}},
\end{equation*}
\begin{equation}\label{eq:getN}
\text{All $k$ measurements}\rightarrow\bY_{(1:k,1)}=\underbrace{\left(\mathcal{O}_{(1,n)}^{(k)}\otimes\mathcal{R}_{(1,n)}^{T(k)}\right)}_{(1:k,n^2)}\underbrace{vec\left(\bN\right)}_{(1:n^2,1)}
\end{equation}

Note that Eqs.\;\eqref{eq:2ndkernelPerDiag},\;\eqref{eq:2ndkernelDiag},\;\eqref{eq:2ndkernel} can be equivalently rewritten as the one linear matrix  equation given in Eq.\;\eqref{eq:getN}. By filling out the above matrix $\left[\mathcal{O}\otimes\mathcal{R}^T\right]$ with the information from $H_{2}(j\omega_1,-j\omega_1)$ and from $H_{2}(j\omega_1,j\omega_1)$ as well, the solution can be improved. Hence, we are able to solve  Eq.\;\eqref{eq:getN} with full rank and identify the matrix $\bN$. All the symmetry properties of the kernels are appropriately used, e.g., conjugate-real symmetry. For $n$ denoting the dimension of the bilinear model and $k$ the number of measurements, we have the following two cases:
\begin{enumerate}
\item $k<n^2$ underdetermined $\rightarrow$ least square (LS) solution (minimizing the 2-norm) as in \cite{KarGA19},
\item $k\geq n^2$ {determined - rank completion} $\rightarrow$ identification of $\bN$,
\end{enumerate}
{\begin{lemma} Let $\Sigma_{b}=(\bA,\bN,\bb,\bc,\bE)$ be a bilinear system of dimension $n$ for which the linear subsystem $\Sigma_{l}=(\bA,\bb,\bc,\bE)$ is \textit{fully controllable and observable}. Then, for $k\geq n^2$ measurements so that $\omega_{1,k},\omega_{2,k}\in\IC$ are distinct pairs, the following holds:
\begin{equation}
rank\underbrace{\left[\begin{array}{c}
\mathcal{O}^{(1)}\otimes\mathcal{R}^{T(1)}\\
\mathcal{O}^{(2)}\otimes\mathcal{R}^{T(2)}\\
\vdots\\
\mathcal{O}^{(k)}\otimes\mathcal{R}^{T(k)}
\end{array}\right]}_{(1:k\geq n^2,n^2)}=n^2
\end{equation}
\end{lemma}
As the above result indicates, one would need at least $n^2$ measurements to identify the matrix $\bN$ corresponding to bilinear system of dimension $n$.}

\subsection{A harmonic separation strategy for the 2nd kernel}\label{sec: separation2}
To identify the $n^{th}$ Volterra kernel, we need an $n$-tone input signal. As we want to identify the 2nd kernel, the input signal needs to be chosen as a double tone Eq.\;\eqref{eq:doubleinput}. The propagating harmonics are: $e^{(j(m_1-m_2)\omega_1+j(m_3-m_4)\omega_2)t}$ or more compactly $e^{(\pm kj\omega_1\pm lj\omega_2)t}$, where $k,l\in\IN$. The aim is to differentiate the $(\omega_1+\omega_2)$ harmonic from the others harmonics.  More precisely, we want the following result to hold:
\begin{equation}\label{eq: phi}
\small
\omega_1+\omega_2\neq k\omega_1+l\omega_2,~\forall(k,l)\in\IZ\times\IZ\setminus{\{1,1\}}.
\end{equation}
Suppose $\omega_2=\phi\omega_1,~\phi\in\IR$. The suitable $\phi$'s where Eq.\;\eqref{eq: phi} holds are:
\begin{equation}\label{eq:eqphi}
\small
\begin{aligned}
\omega_1+\phi\omega_1&=k\omega_1+l\phi\omega_1\Rightarrow1+\phi=k+l\phi\Rightarrow\phi=\frac{k-1}{1-l},~k,l\in\IZ\setminus\{1\}.
\end{aligned}
\end{equation}

By choosing $\phi$ so that the equality in Eq.\;\eqref{eq:eqphi} doesn't hold, with harmonic mixing index $m=k+l$ (i.e., for 2nd stage approximation $m=2$), it makes the harmonic $(\omega_1+\omega_2)$ uniquely defined in the frequency spectrum up to the $m^{th}$ harmonic. 

\begin{figure}
\centering
\includegraphics[scale=0.8]{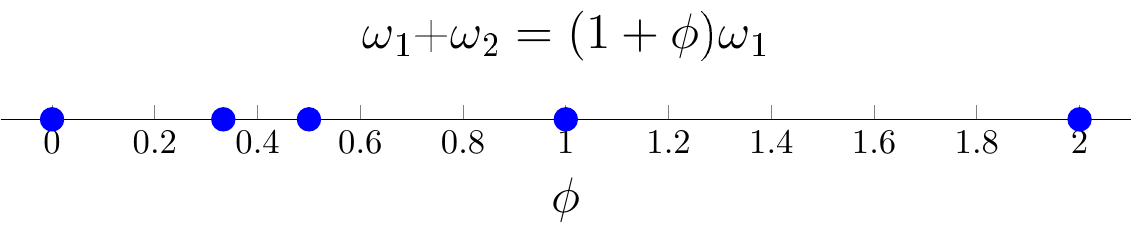}
\caption{This figure shows the constrains of $\phi$ (e.g., $\phi=0,1/3,1/2,1,2,3,\ldots$, etc.). By choosing $\phi$'s within the blue dots, we construct frequency bandwidths with an unique $(\omega_1+\omega_2)$.}
\label{fig:Phi1dviolations}
\end{figure}

To visualize this feature, we choose $\omega_1=1$, and $\omega_2=\omega_1\phi=\phi$, for harmonic mixing index $m=2$. Then, the constraints of $\phi$ are depicted in Fig.\;\ref{fig:Phi1dviolations} with blue dots.

\begin{figure}
\centering
\includegraphics[scale=0.45]{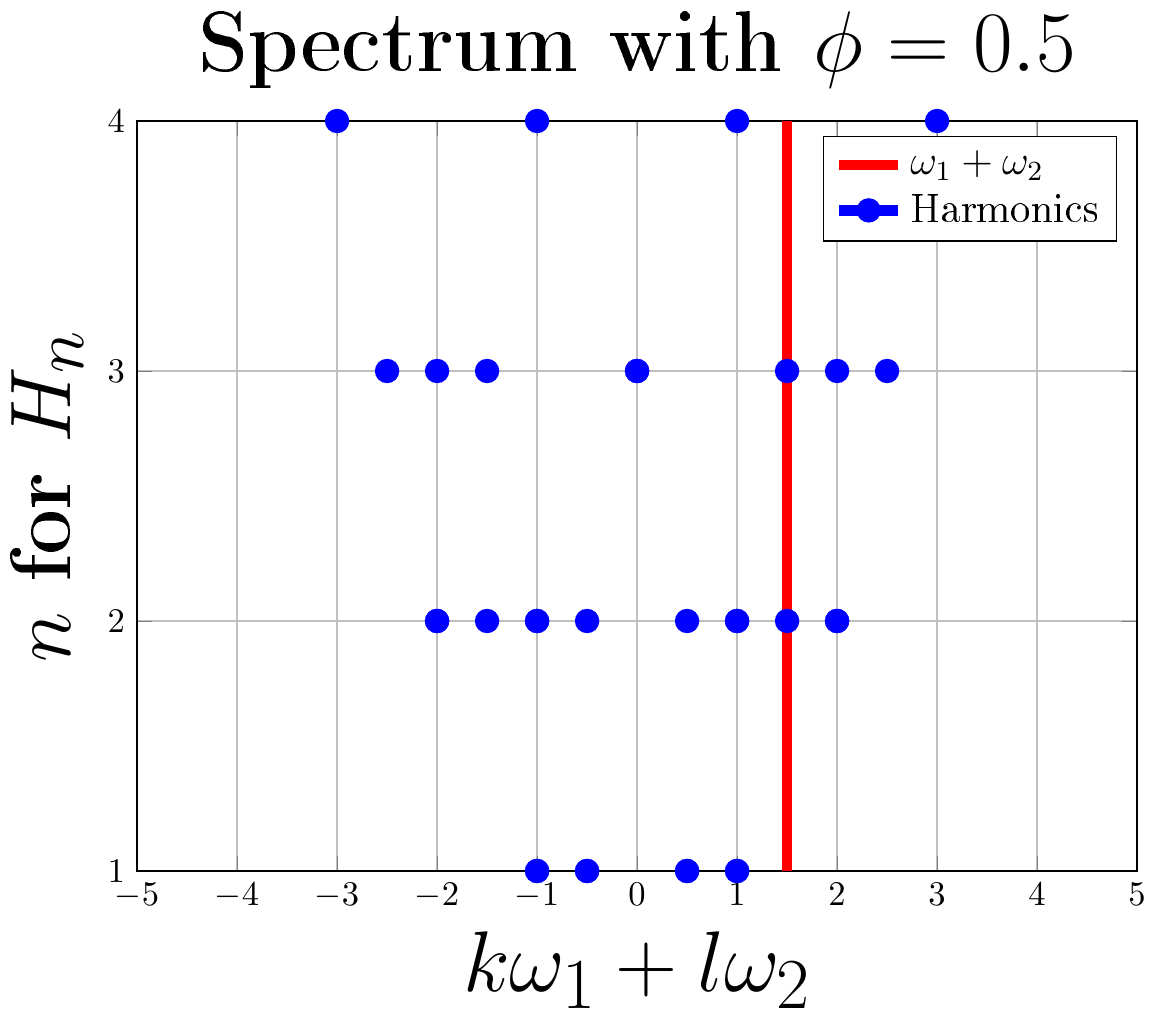}
\includegraphics[scale=0.45]{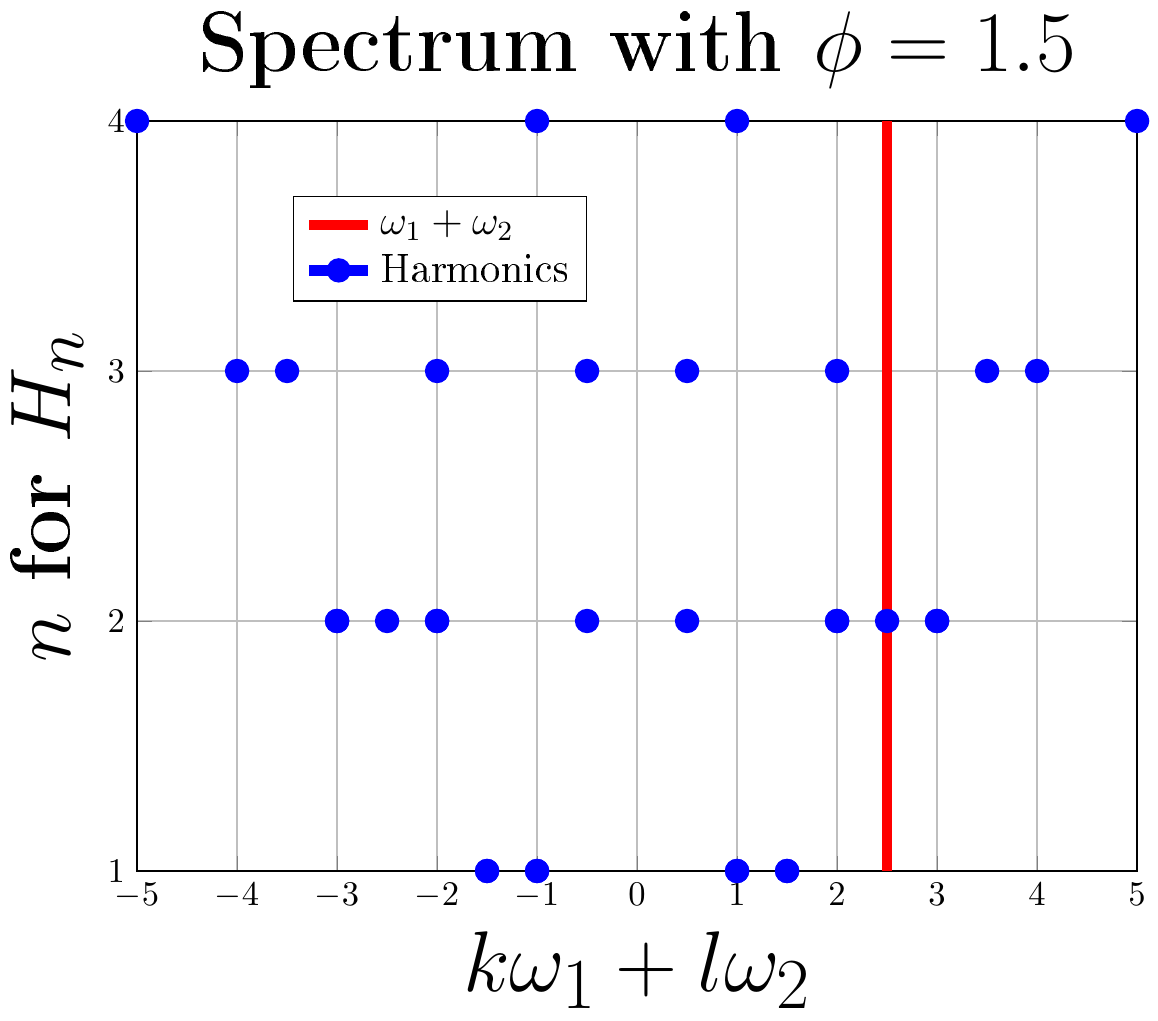}
\caption{Left pane: Overlapping harmonics with invalid $\phi=0.5$. Right pane: Uniquely defined harmonic at $(\omega_1+\omega_2)$ with valid $\phi=1.5$. Here, it holds $(n=k+l)$.}
\label{fig:violations}
\end{figure}

Next, in Fig.\;\ref{fig:violations} and on the left pane, one $\phi$ constraint that occurs commensurate harmonics is depicted,  and on the right pane, one unique harmonic construction at $(\omega_1+\omega_2)$ is presented.

The next result allows us to construct sweeping frequency schemes to get enough measurements for the $H_{2}(j\omega_1,j\omega_2)$. So, for every $\omega_1>0$ the following should hold:
\begin{equation}
\small
\omega_2\in\left(\phi_{i-1}\omega_1,\phi_{i}\omega_1\right),~i=1,\ldots
\end{equation}
where $\phi_i$ are the constrains (see Fig. 4 blue dots).
\begin{remark}
Note that in the proposed framework, the separation of the $\omega_1+\omega_2$ harmonic is forced only under a specific mixing order $m$. We do not offer any general solution of the harmonic separation problem for multi-tone input, although techniques have been introduced such as in \cite{Boyd1983MeasuringVK}. There, it was also stated that the solution of the full separation of harmonics is in general, not possible.
\end{remark}
\subsection{The Loewner-Volterra algorithm for time-domain bilinear identification}
We start with a set of single-tone inputs $u(t)=\alpha_1\cos(\omega_i t),~i=1,...,n$, with $\alpha_1<1$. For those $n$ measurements we can estimate the linear kernel $H_{1}(j\omega_i)$, the $H_{2}(j\omega_{i},j\omega_{i})$ and the $H_{2}(j\omega_i,-j\omega_i)$ by simply measuring the first harmonic as $\bY_{1}$, the second harmonic as $\bY_{2}$, and the DC term as $\bY_{0}$, from the frequency spectrum as in Fig.\;\ref{fig:Spectrum}. To improve the accuracy of the estimations for the aforementioned kernels, we could further upgrade to an $\ell$-stage approximation by varying the amplitude $\alpha_{\ell}$ as explained in section \ref{sec: separation}. This approach is necessary whenever higher harmonics are considered to be numerically non-zero, hence meaningful. The reason for this is that the first harmonic is hence corrupted by noise introduced by the term $H_{3}^{2,1}$ which appears on the second row of matrix $\bP_\omega$ in Eq.\;\eqref{eq:system_YMP}.

Since the LF reveals the underlying order of the system denoted with $r$, the value of $n$ should be at least equal to $2r$. Then, we can take the decision on what will be the order $r$ of the reduced system by analysing the singular value decay. Up to the previous step, we have identified the linear part with the LF, and, we have filled the LS problem Eq.\;\eqref{eq:getN} with measurements from the diagonal of the second kernel and from the the perpendicular to the diagonal axis $(\omega_1,-\omega_1)$. Those measurements contribute to the LS problem, but with an underdetermined (rank deficient) LS problem.  

We need more measurements of $H_{2}$ to reach the full rank $(r^2)$ solution that will lead to the identification of $\bN$. So, we proceed by measuring the $H_{2}$ off the diagonal with a double-tone input as $u(t)=\alpha_{\ell}\cos(\omega_{1,k} t)+\beta_{\ell}\cos(\omega_{2,k} t)$, for a set of frequency pairs $(\omega_1,\omega_2)$ up to $r^2$. The harmonic separation problem for the frequency $(\omega_1,\omega_2)$ appears now. To deal with this problem, we follow the solution proposed in section \ref{sec: separation2} (up to a mixing degree). Last, we solve the real\footnote{Enforcing real-valued models has been discussed in \cite{tutorial,KGAHandbook}; here, we follow the same approach.} full rank LS problem described in Eq.\;\eqref{eq:getN} by using all the symmetric properties of these kernels (i.e., real symmetry, conjugate symmetry and the fact that $H_{2}(j\omega_1,j\omega_2)=H_{2}(j\omega_2,j\omega_1)$). An algorithm that summarizes the above procedure is presented below.

\begin{programcode}{(Algorithm 2) A Loewner-Volterra algorithm for bilinear identification/reduction from time-domain data.}
\vspace{-5mm}
\begin{enumerate}
\item[] \textbf{Input/{Data acquisition}:} Use as control input the signals: $u(t)=\alpha_{\ell}\cos(\omega_{1,k} t)+\beta_{\ell}\cos(\omega_{2,k} t),~t\geq 0$, by sweeping the small amplitudes $(<1)$ and a particular range of frequencies.\\
\item[] \textbf{Output:} A reduced bilinear system of dimension-$r$: $\Sigma_{r}:\left(\bA_{r},\bN_{r},\bb_{r},\bc_{r},\bE_{r}\right)$ \\[-2mm]
\item Apply one-tone input $u(t)$ with $\beta_{\ell}=0$, $\omega_{1,k}$ for $k=1,\ldots,n$ and collect the snapshots $y(t)$ in steady state. Improve the measurements by solving system Eq.\;\eqref{eq:amplitudeImprovement} by sweeping the amplitude for each frequency.\\
\item Apply Fourier transform for each frequency and collect the following $k$ measurements:
\begin{itemize}
\item DC term: $~Y_{O}(0\cdot j\omega_{1,k})$,
\item 1st harmonic: $Y_{I}(1\cdot j\omega_{1,k})$,
\item 2nd harmonic: $Y_{II}(2\cdot j\omega_{1,k})$.
\end{itemize}
\item Apply the linear LF, see Algorithm 1 in \cite{KGAHandbook} by using the measurements (e.g., $H_{1}(j\omega_{1,k})=2Y_{I}(j\omega_{1,k})/\alpha_{1,\ell}$ for the 2nd stage approximation) and get the the order $r$ reduced linear model. If the 2nd harmonic and higher harmonics are equal with zero, the underlying system is linear so $\bN$ is the zero matrix $\bf{0_{r\times r}}$ as well. Otherwise the system is nonlinear and computing the bilinear matrix $\bN$ will improve the accuracy.\\
\item Apply the 2-tone input $u(t)=\alpha_{\ell}\cos(\omega_{1,k} t)+\beta_{\ell}\cos(\omega_{2,k} t)$ to get enough measurements $(\leq r^2)$ to produce a full rank LS problem. Measure the $(\omega_1+\omega_2)$ harmonic as explained in section \ref{sec: separation2} and get the estimations for the 2nd kernel as: $H_{2}(j\omega_{1,k},j\omega_{2,k})=2Y_{II}(j\omega_{1,k},j\omega_{2,k})/(\alpha_{\ell}\beta_{\ell})$.\\
\item Solve the full rank least square problem described in Eq.\;\eqref{eq:getN} and compute the real-valued bilinear matrix $\bN$. {When the inversion is not exact due to numerical issues, the least squares solution is obtained with a thresholding SVD.}
\end{enumerate}
\end{programcode}
\subsection{Computational effort of the proposed method}
In this section we discuss the computational effort of the proposed method by analyzing each step. We comment on the applicability of large scale problems and the relation with real-world scenarios.

Simulation of processes with harmonic inputs constitutes a classical technique which is applied in many engineering applications; data acquisition in the time domain is a common procedure. Nevertheless, using advanced electronic devices such as vector network analyzers (VNAs), frequency domain data can also be obtained (directly). The Loewner framework applied in the case where frequency domain data that are obtained from VNAs offers an excellent identification and reduction tool in the linear case (with many applications in electrical, mechanical or civil engineering). In the context of the current paper, we deal with time-domain data for a special class of nonlinear problems.

For the purpose of identifying and reducing bilinear systems from time-domain measurements, the most expensive procedure is collecting the data. This is done by simulating time-domain models with Euler's method (bilinear models such as the ones approximating Burgers' equation). Nevertheless, the heavy computational cost of simulating large dimensional systems in time domain  could be alleviated using parallel processing (e.g., for multiple computational clusters). The process of estimating transfer functions values by computing the Fourier transform hence remains robust. In addition, the LF can adaptively detect the decay of the singular values and hence the procedure can be terminated for a specific reduced order $r\ll n$.

In the beginning, a linear system of reduced dimension $r$ is fitted using the LF. For the rest of the proposed algorithm, note that we will use the lower dimension $r$ to our advantage, and hence, the method remains robust. The next step is to compute the matrix $\bN$ that characterizes the nonlinearity of bilinear systems. As the fitted linear system is of dimension $r$, we hence need to detect exactly $r^2$ unknowns (the entries of matrix $\bN$).As presented in section 4.3, this boils down to solving a full-rank LS problem that can be easily dealt with.

The aim of the newly-proposed method is to accurately train bilinear models from time-domain data.  We offer a first step approach towards complete identification of such systems within the Volterra series approximation approach. In many cases, large-scale systems are sparse (due to spatial domain semi-discretization) and hence, reduction techniques can be applied. The new method deals with the inherent redundancies through the linear subsystem (compression by means of SVD). Afterwards, it updates the nonlinear behavior by introducing an appropriate low-dimensional bilinear matrix that improves the overall approximation. Note also that the new method relies on the \textit{controllability/observability} of the fitted linear system. Additionally, noise values (due to FFT) up to a particular threshold can be handled as presented in section 5; further analysis on noise-related issues is left for future research.

\section{Numerical examples}
\begin{example}[Identifying a low-order bilinear toy example] The aim of this experiment is to identify a simple bilinear model from time-domain measurements. Consider the following controllable/observable bilinear model Eq.\;\eqref{eq:bilinear} of dimension-$2$ with a \textit{non-symmetric} matrix $\bN$, zero initial condition and matrices as:
\begin{equation}
\small
\bE=\left[\begin{array}{cc}
1 & 0\\
0 & 1
\end{array}\right],~\bA=\left[\begin{array}{cc} -1 & -10\\ 10 & -1 \end{array}\right],~\bN=\left[\begin{array}{cc} 1 & -2\\ 3 & -4 \end{array}\right],~\bB=\left[\begin{array}{c} 1\\ 1 \end{array}\right],~
\bC=\left[\begin{array}{cc} 1 & 1 \end{array}\right].
\end{equation} 
We simulate the system in the time domain with a single sinusoidal input $u(t)=0.01\cos(2\pi\omega t)$, by varying the frequency $\omega\in [0.5,1,1.5,2]$, with time step $dt=1e-4$. Next, the 2nd-stage approximation results for the linear kernel $\tilde{H}_{1}$ in comparison with the theoretical values of $H_{1}$ are presented in Table\;\ref{tab:H1}.
\begin{table}[!h]
\caption{Measurements of the first (linear) kernel}
\label{tab:H1}       
\begin{tabular}{p{3cm}p{3cm}p{3cm}}
\hline\noalign{\smallskip}
Frequency $\omega$(Hz) & $\hat{H}_1(j\omega)$-2nd stage & $H_{1}(j\omega)$-theoretical \\
\noalign{\smallskip}\svhline\noalign{\smallskip}
0.5 & $+0.026606+0.067106{}\mathrm{i}$ & $+0.026574+0.067115{}\mathrm{i}$\\
1 & $+0.071503+0.189600{}\mathrm{i}$  & $+0.071258+0.189700{}\mathrm{i}$\\
1.5 & $+0.752720+0.377300{}\mathrm{i}$  & $+0.754030+0.380870{}\mathrm{i}$\\
2 & $+0.134070-0.381970{}\mathrm{i}$  & $+0.133780-0.382520 {}\mathrm{i}$\\
\noalign{\smallskip}\hline\noalign{\smallskip}
\end{tabular}\\
$^a$ With 2nd-stage approximation $\tilde{H}_{1}(j\omega)= 2Y_{1}(j\omega)/\alpha$.
\end{table}

With the estimations of the linear transfer function and by using the (LF) as the data-driven identification and reduction tool for linear systems, we identify the linear system $(\tilde{\bA},\tilde{\bb},\tilde{\bc},\tilde{\bE})$. We stopped at the 4th measurement due to the fact that the underlying system is of second order (McMillan degree 2). Otherwise, more measurements will be needed to have a sufficient decay of the singular values as in Fig.\;\ref{fig:loewnerDecay}. The singular values decay offers a choice of reduction. As long as the simulation of the system is done, with time step $dt=1e-4$, the singular values with magnitude below that threshold are neglected. 
\begin{figure}[!h]
\centering
\includegraphics[scale=0.4]{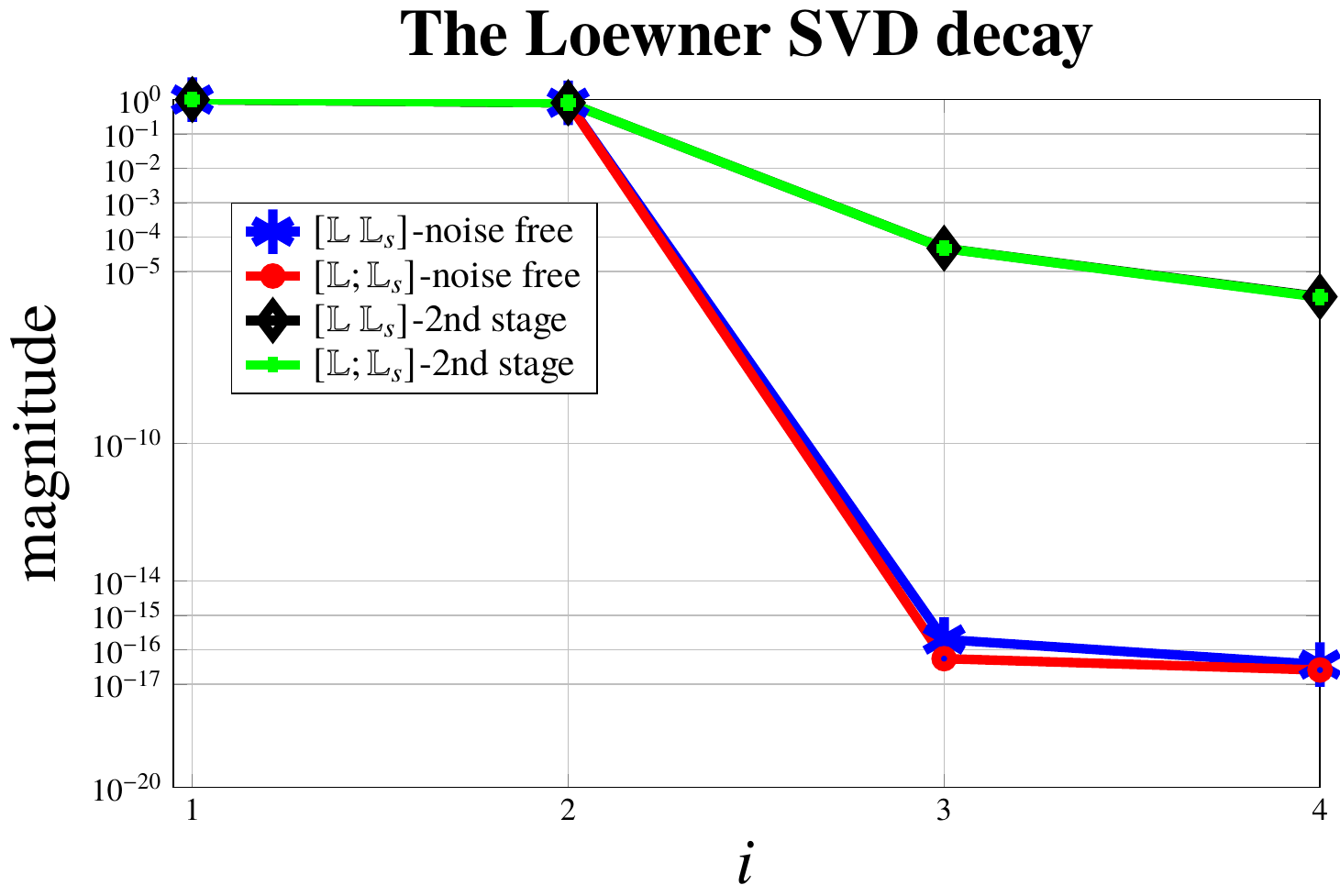}
\caption{The singular value decay of the (LF) as a fundamental characterization of the McMillan degree of the underlying linear system. Here, a truncation scheme of order $r=2$ is recommended where the 2nd stage approximation gave $\sigma_{3}/\sigma_{1}=4.721\cdot 10^{-5}$, while for the noise free case the third singular values has reached the machine precision.}
\label{fig:loewnerDecay}
\end{figure}

Construction of the linear system with order $r=2$ by using the theoretical noise-free measurements (subscript "t") appears next:
\begin{equation}
\small
\tilde{\bA}_t=\left[\begin{array}{cc} -1.4513 & -8.8181\\ 11.363 & -0.54868 \end{array}\right],~\tilde{\bB}_t=\left[\begin{array}{c} -0.92979\\ 1.3967 \end{array}\right],~
\tilde{\bC}_t=\left[\begin{array}{cc} -0.76857 & 0.9203 \end{array}\right],
\end{equation}
while by using the measured data with 2nd-stage approximation results to the following:
\begin{equation}
\small
\tilde{\bA}=\left[\begin{array}{cc} -1.458 & -8.8137\\ 11.367 & -0.55162 \end{array}\right],~\tilde{\bB}=\left[\begin{array}{c} -0.9342\\ 1.4 \end{array}\right],~
\tilde{\bC}=\left[\begin{array}{cc} -0.7675 & 0.91611 \end{array}\right].
\end{equation}
\begin{important}{Identified linear dynamics}
Even if the coordinate system is different, one crucial qualitative result is to compute the poles and zeros of the linear transfer function. For the identified system with the theoretical measurements (noise free), the poles and zeros are exactly as the original: $\tilde{p}_t=-1\pm 10{}\mathrm{i}$ and the zero is: $\tilde{z}_t=-1$ while for the 2nd stage approximation to the linear system, the corresponding results are: $\tilde{p}=-1.0048\pm 9.9989{}\mathrm{i},~\tilde{z}=-1.0042$.
\end{important}
At this point, we have recovered the linear part of the bilinear system up to an accuracy due to the truncation of Volterra series. The inexact simulations of the continuous system which are done with a finite time step $dt=1e-4$, and the Fourier accuracy led to quite accurate results with a perturbation of the order $\sim O(1e-3)$ by comparing the theoretical poles and zeros. We proceed by collecting the measurements of the 2nd
 kernel. The following Table\;\ref{tab:2}, contains measurements of the 2nd kernel with 1-tone input.
\begin{table}[!h]
\caption{Measurements of the $H_{2}$ on the diagonal and perpendicular to the diagonal.}
\label{tab:2}       
\begin{tabular}{p{1.3cm}p{3cm}p{3cm}p{2cm}p{2cm}}
\hline\noalign{\smallskip}
Frequency $\omega$(Hz) & $\tilde{H}_2(j\omega,j\omega)$ & $H_{2}(j\omega,j\omega)$-theoretical &$\tilde{H}_2(j\omega,-j\omega)$& $H_2(j\omega,-j\omega)$-theoretical \\
\noalign{\smallskip}\svhline\noalign{\smallskip}
0.5 & $+0.026440-0.124490{}\mathrm{i}$  & $+0.026570-0.124440{}\mathrm{i}$ & $+0.032190$ & $+0.032177$ \\
1 & $-0.184590+0.298430{}\mathrm{i}$  & $-0.184510+0.298910{}\mathrm{i}$ & $+0.045648$ & $+0.045641$ \\
1.5 & $+0.178080+0.305840{}\mathrm{i}$  & $+0.178160+0.307170{}\mathrm{i}$ & $+0.063936$ & $+0.064350$ \\
2 & $+0.062642-0.054219{}\mathrm{i}$  & $+0.062588-0.054423{}\mathrm{i}$ & $-0.044927$ & $-0.044998$ \\
\noalign{\smallskip}\hline\noalign{\smallskip}
\end{tabular}\\
$^b$ The estimation of the 2nd kernel with 2nd stage approximation on the diagonal as $\tilde{H}_{2}(j\omega,j\omega)=2Y_{2}(j\omega,j\omega)/\alpha^2$ and perpendicular of the diagonal as $\tilde{H}_{2}(j\omega,-j\omega)=2Y_{2}(j\omega,-j\omega)/\alpha^2$.
\end{table}

We can get $\bN$ by solving the least square problem by just minimizing the 2-norm as in \cite{KarGA19}. This result was not towards the identification of the matrix $\bN$ and here is the  new approach working towards the identification of bilinear systems.

\begin{question}{Can we identify the matrix $\bN$?}
The improvement here relies on the rank deficiency problem is produced by solving the least square problem without taking under consideration measurements off the diagonal of the 2nd kernel $H_{2}$ along $(\omega_1,\omega_2)$ with ($\omega_1\ne\omega_2$). By filling the LS problem Eq.\;\eqref{eq:getN} with a full rank, the answer is affirmative.
\end{question} 
Back to our introductory example, the rank of the least squares problem is less than $r^2=4$. So, we need to increase the rank. We take measurements ($\leq 4$) off the diagonal from the 2nd kernel by using the 2-tone input. Tab.\;\ref{tab:4} includes the theoretical and measured results.
\begin{table}[!h]
\caption{Measurements of the 2nd kernel (off the diagonal)}\label{tab:4}   
\begin{tabular}{p{4cm}p{3cm}p{3cm}}
\hline\noalign{\smallskip}
Frequencies $(\omega_1,\omega_2)$(Hz) & $\tilde{H}_{2}(j\omega_1,j\omega_2)$ & $H_{2}(j\omega_1,j\omega_2)$ \\
\noalign{\smallskip}\svhline\noalign{\smallskip}
$(0.2,0.3)$ & $+0.030440-0.039259{}\mathrm{i}$  & $+0.030429-0.039237{}\mathrm{i}$ \\
$(0.2,0.6)$ & $+0.031002-0.080364{}\mathrm{i}$  & $+0.031037-0.080315{}\mathrm{i}$ \\
$(0.4,0.3)$ & $+0.030948-0.062869{}\mathrm{i}$  & $+0.030961-0.062835{}\mathrm{i}$ \\
$(0.4,0.6)$ & $+0.026417-0.125320{}\mathrm{i}$  & $+0.026554-0.125260{}\mathrm{i}$ \\
\noalign{\smallskip}\hline\noalign{\smallskip}
\end{tabular}\\
$^c$ The estimation of the 2nd kernel as $\tilde{H}_{2}(j\omega_1,j\omega_2)=Y_{2}(j\omega_1,j\omega_2)/(\alpha_1\alpha_2)$. Here we use $\phi=1.5$, to avoid the harmonic overlapping as explained in section \ref{sec: separation2}.
\end{table}

The full rank LS solution gave for the theoretical noise free case and for the 2nd stage approximation the following results respectively:
\begin{equation}
\small
\tilde{\bN}_t=\left[\begin{array}{cc} -4.1542 & -2.0998\\ 3.236 & 1.1542 \end{array}\right],~\tilde{\bN}=\left[\begin{array}{cc} -4.1557 & -2.1084\\ 3.2284 & 1.1513 \end{array}\right]
\end{equation}
\vspace{-5mm}
\begin{important}{Coordinate transformation}
By transforming all the matrices to the same coordinate system as in \cite{Juang}, we conclude:
\begin{itemize}
\item \textbf{Noise free case - exact identification}\\
\begin{equation}
\small
\breve{\bA}_t=\left[\begin{array}{cc} -1.0 & -10.0\\ 10.0 & -1.0 \end{array}\right],~\breve{\bN}_t=\left[\begin{array}{cc} 1.0 & -2.0\\ 3.0 & -4.0 \end{array}\right],~\breve{\bB}_t=\left[\begin{array}{c} 1.0\\ 1.0 \end{array}\right],~
\breve{\bC}_t=\left[\begin{array}{c} 1.0 \\ 1.0 \end{array}\right]^T.
\end{equation}
\item \textbf{Simulated case - approximated identification}\\
\begin{equation}
\small
\breve{\bA}=\left[\begin{array}{cc} -1.0037 & -9.9941\\ 10.004 & -1.0059 \end{array}\right],~\breve{\bN}=\left[\begin{array}{cc} 0.99525 & -1.997\\ 3.006 & -3.9997 \end{array}\right],~\breve{\bB}=\left[\begin{array}{c} 0.99925\\ 1.0003 \end{array}\right],~
\breve{\bC}=\left[\begin{array}{c} 1.0 \\ 1.0 \end{array}\right]^T.
\end{equation}
\end{itemize}
\end{important}
Next, in Fig.\;\ref{fig:kernels}, evaluation results for the linear and the second order  generalized transfer function are presented:
\begin{figure}[!h]
\centering
\includegraphics[scale=0.25]{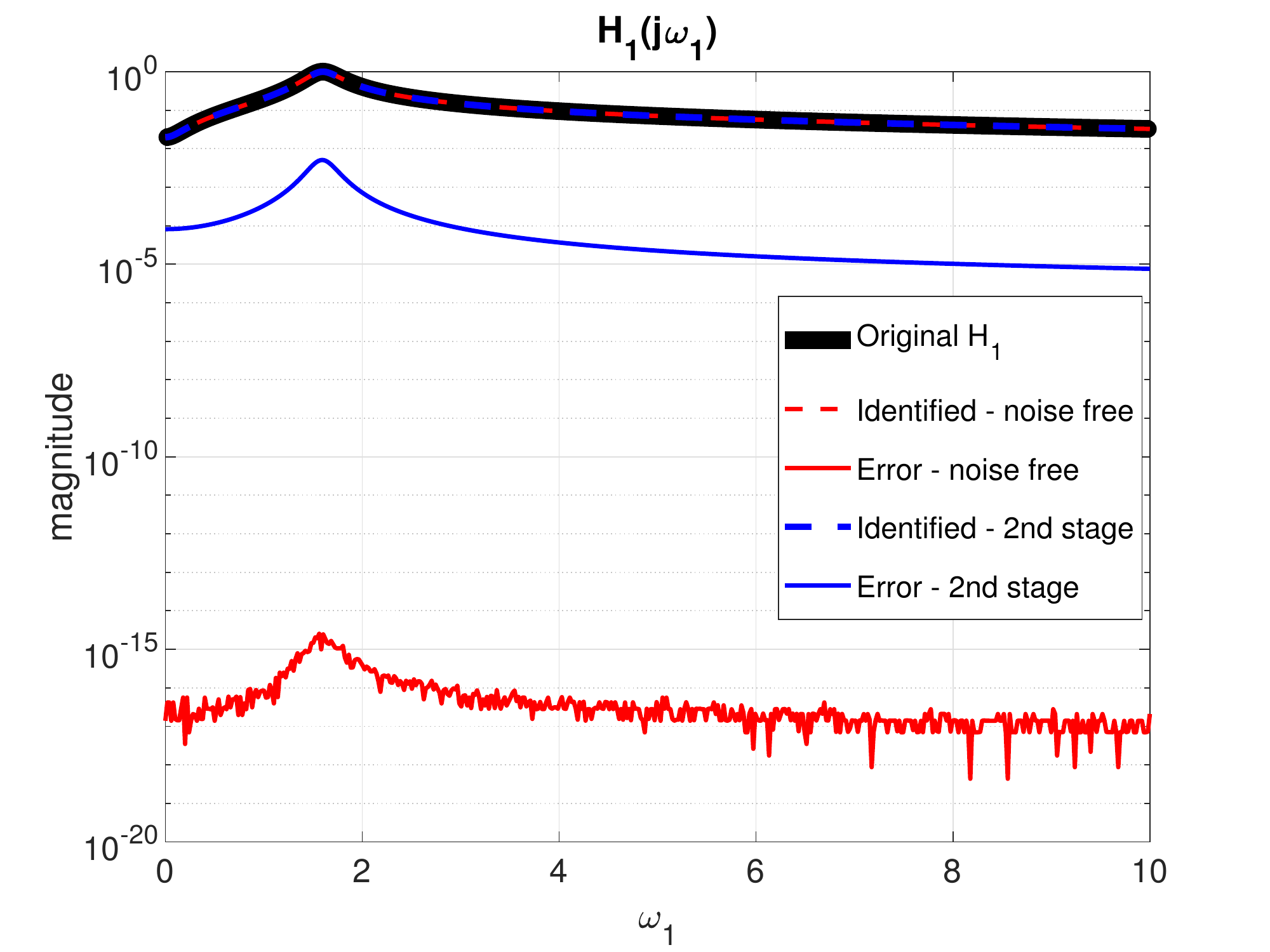}
\includegraphics[scale=0.25]{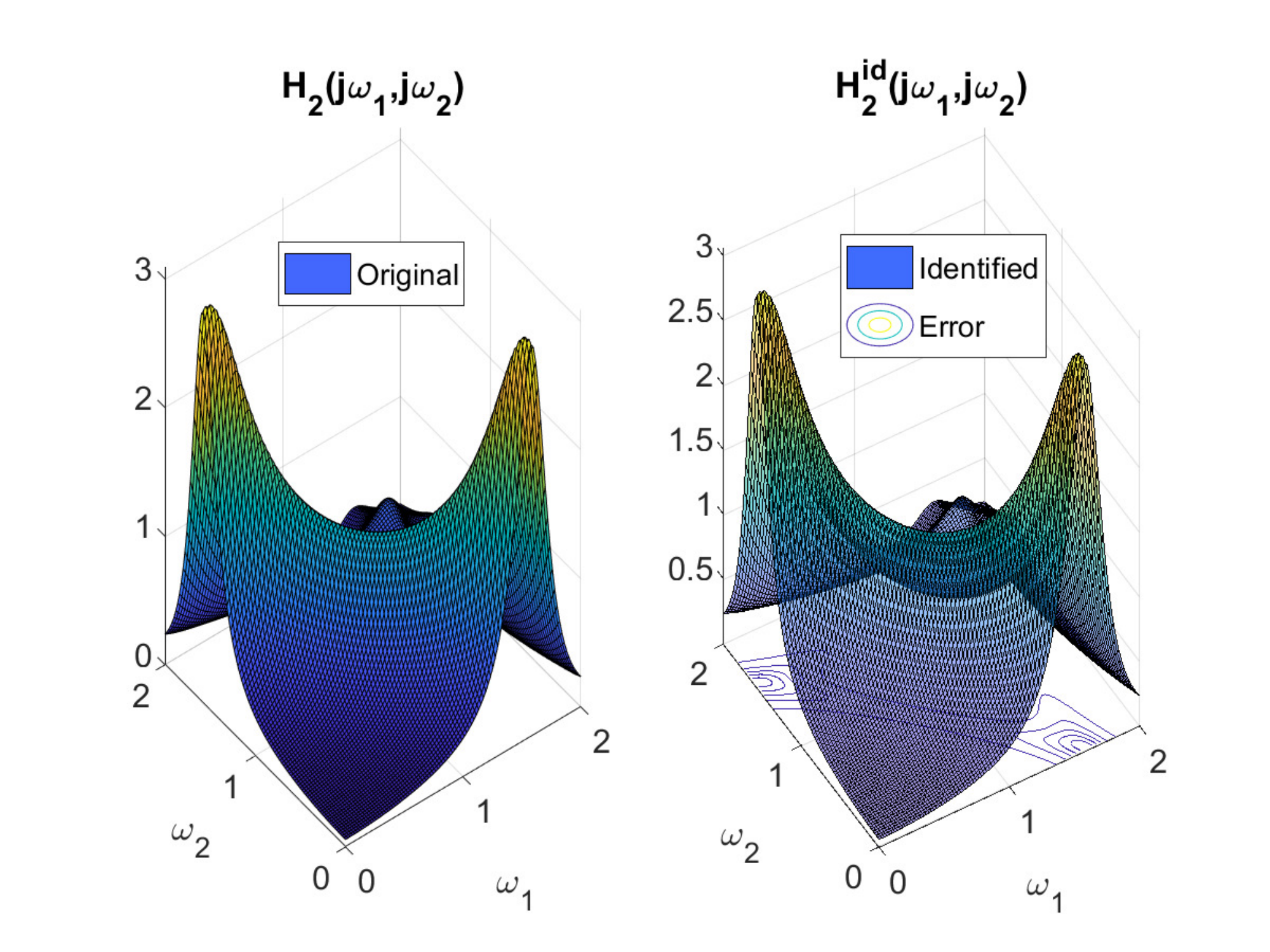}
\caption{The identified 1st and 2nd kernel with 2nd-stage approximation in comparison with the theoretical kernels.}
\label{fig:kernels}
\end{figure}

Finally, time-domain simulations for each system performed in Fig.\;\ref{fig:evalut} with a larger amplitude than the probing one.
\begin{figure}[!h]
\centering
\includegraphics[scale=0.5]{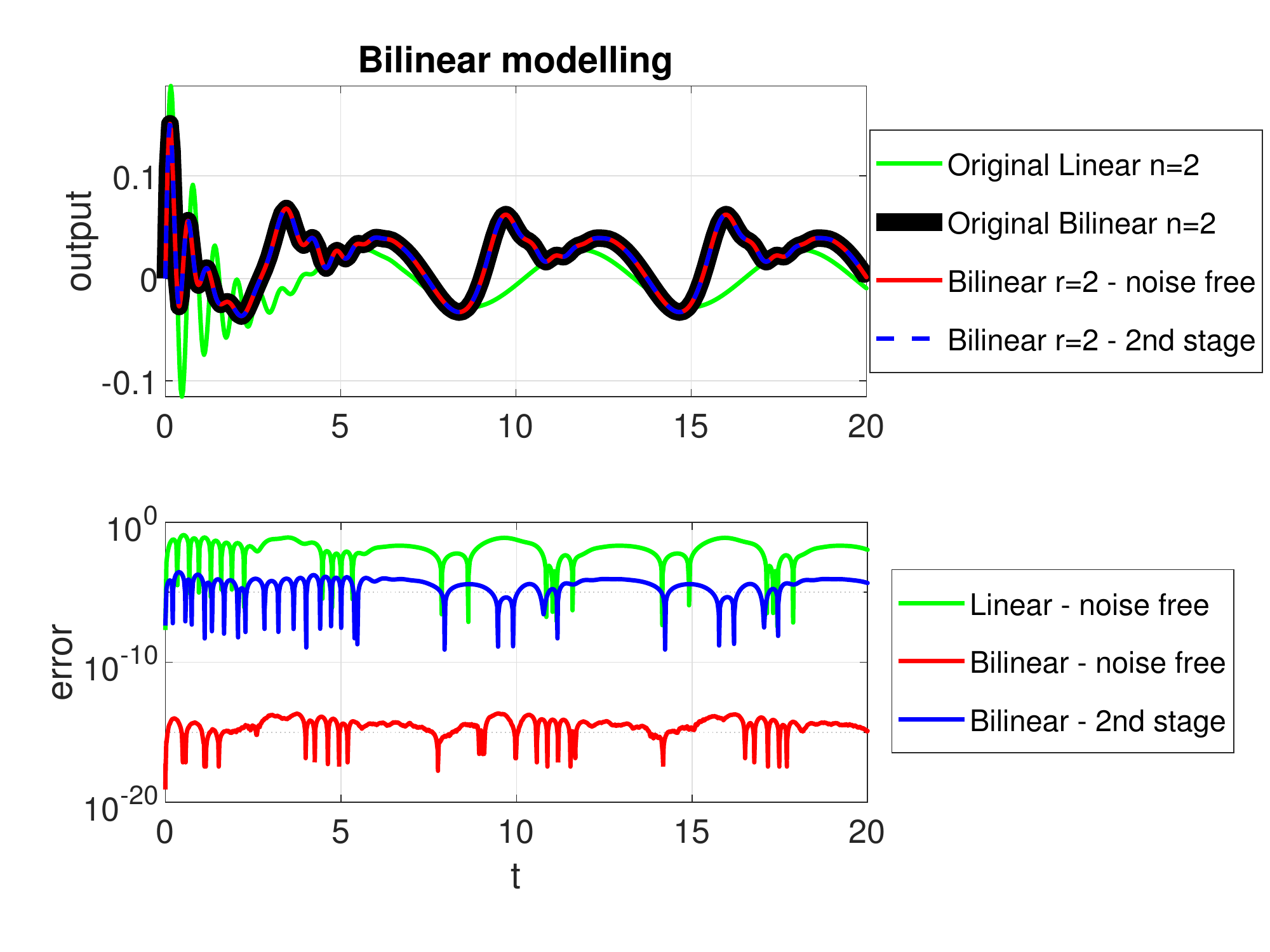}
\caption{The evaluation of the models with order $r=2$ performed with input as $u(t)=\cos(t),~t\in [0,20]$. The noise free case has reached machine precision.}
\label{fig:evalut}
\end{figure}
\end{example}
\begin{example}
{\textbf{Time domain reduction of the Burgers' Equation.}}
This example illustrates the bilinear modelling and reduction concepts proposed in \cite{AGIbil} for the viscous Burgers' equation from time-domain simulations. We simulate the system with $40$ measurements as $\omega_{k}=j2\pi[0.1,0.2,\ldots,4]$. We present the corresponding results with initial system dimension $n=420$ reduced by the proposed method to order $r=2$ with the first neglected singular value to be $\sigma_{3}/\sigma_{1}=4.6255\cdot 10^{-4}$. As the order was chosen $r=2$, the reduced bilinear matrix $\tilde{\bN}$ was introduced by using the following measurements as $\omega_{1}=j2\pi[0.2,0.4]$ and $\omega_{2}=j2\pi[0.3,0.6]$. In Fig.\;\ref{fig:burgersKernels}, evaluation results are presented.
\begin{figure}[!h]
\centering
\includegraphics[scale=0.25]{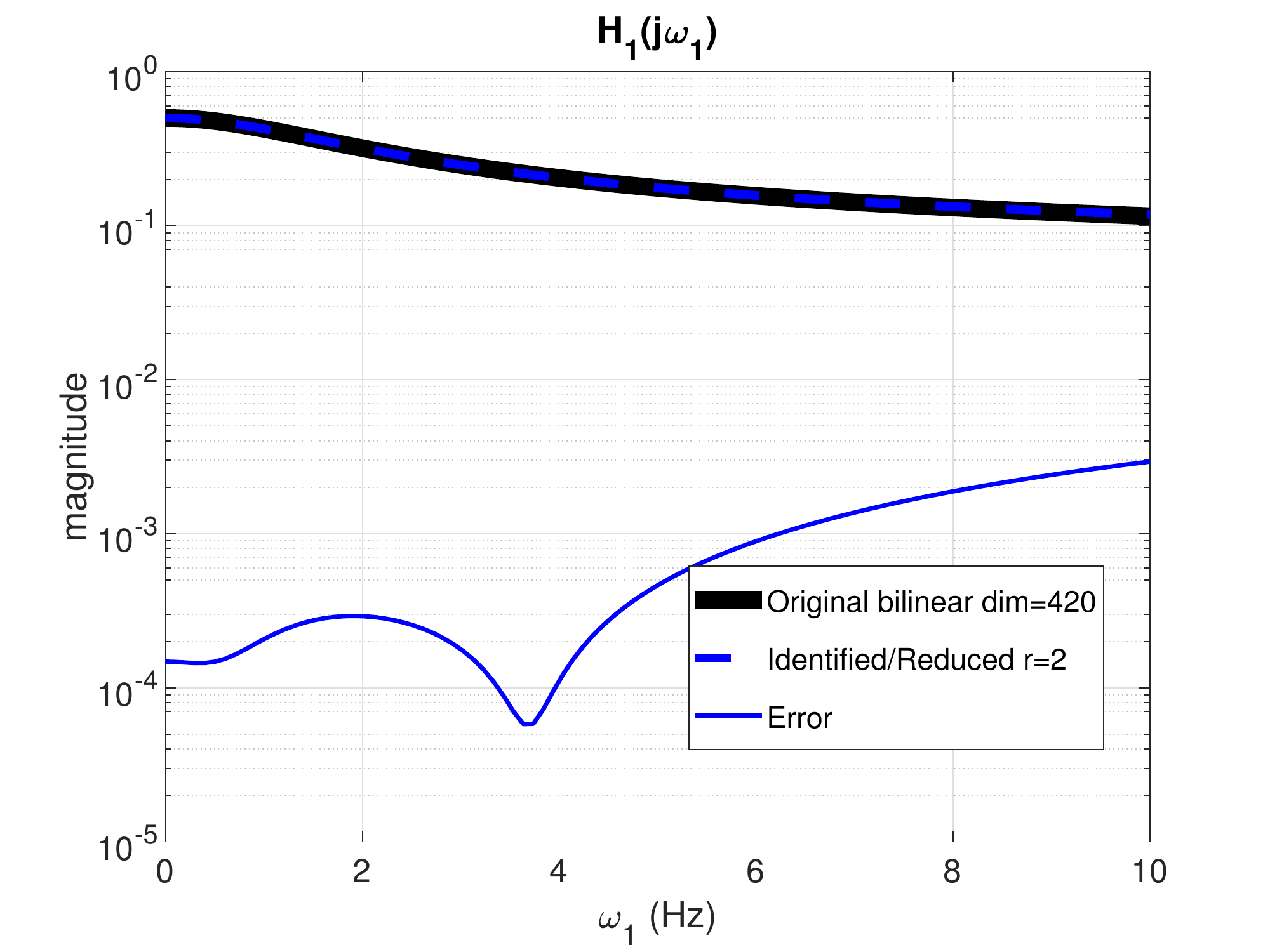}
\includegraphics[scale=0.25]{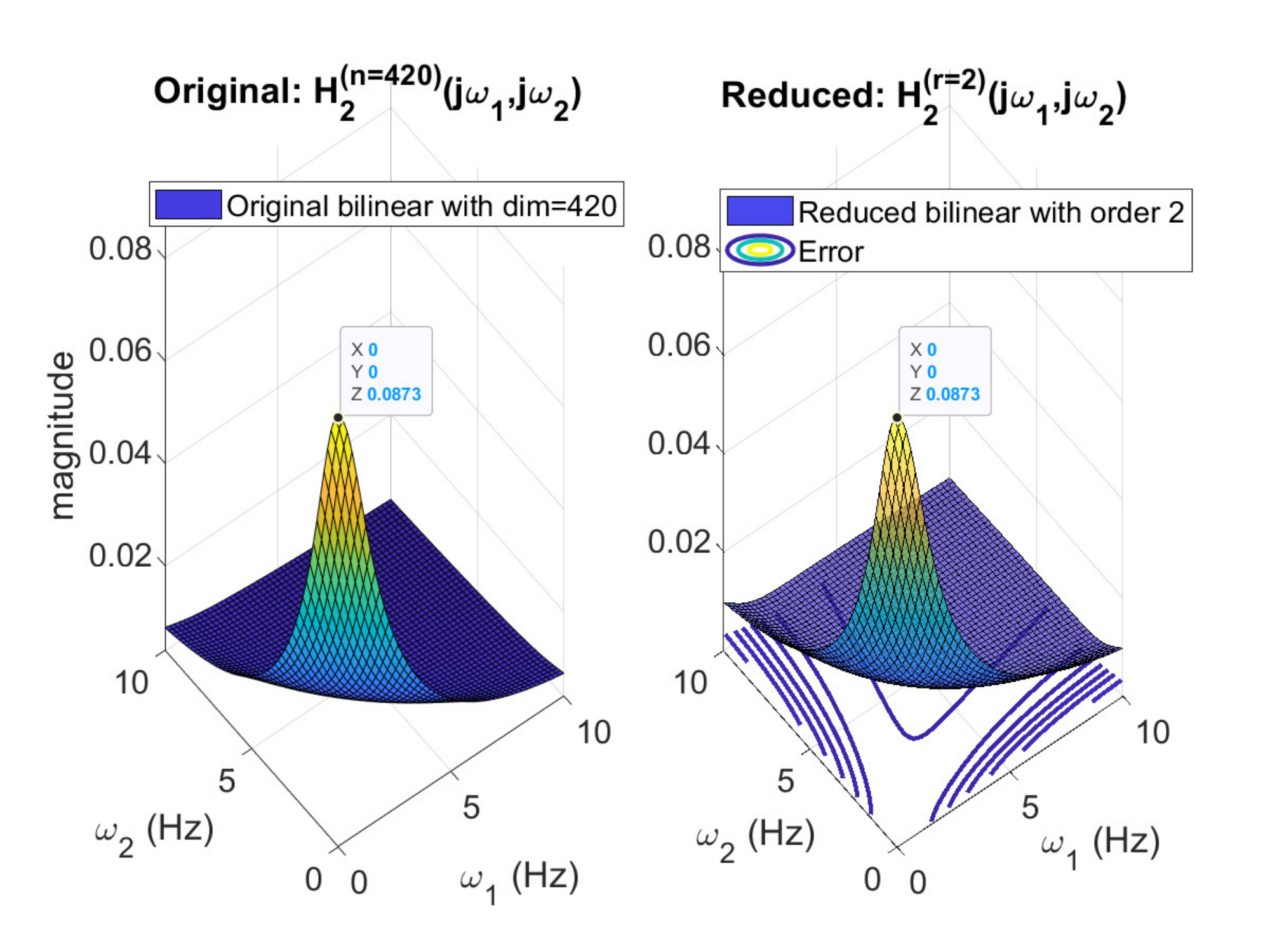}
\caption{The 1st and the 2nd kernel evaluations in comparison with the originals.}
\label{fig:burgersKernels}
\end{figure}
Lastly, in Fig.\;\ref{fig:burgersSimulation}, a time-domain simulation reveals that the proposed method can improve the accuracy by fitting a nonlinear model. Table 4 contains approximation results both in the frequency and, also in the time domain. For the example presented (dimension reduction from $n=420$ to $r=2$), we offer a comparison of the newly-proposed method (Time-LoewBil) with another method, i.e., the frequency-domain bilinear Loewner framework introduced in \cite{AGIbil} (Freq-LoewBil). The common frequency grid was selected as described above while the sampling values of the tranfser functions (in the frequency domain) were corrupted with white-noise. The noise magnitude of the latter was selected to match the noise values introduced by performing time-domain simulations with a time step of $dt=1e-4$.

\begin{figure}[!h]
\centering
\includegraphics[scale=0.5]{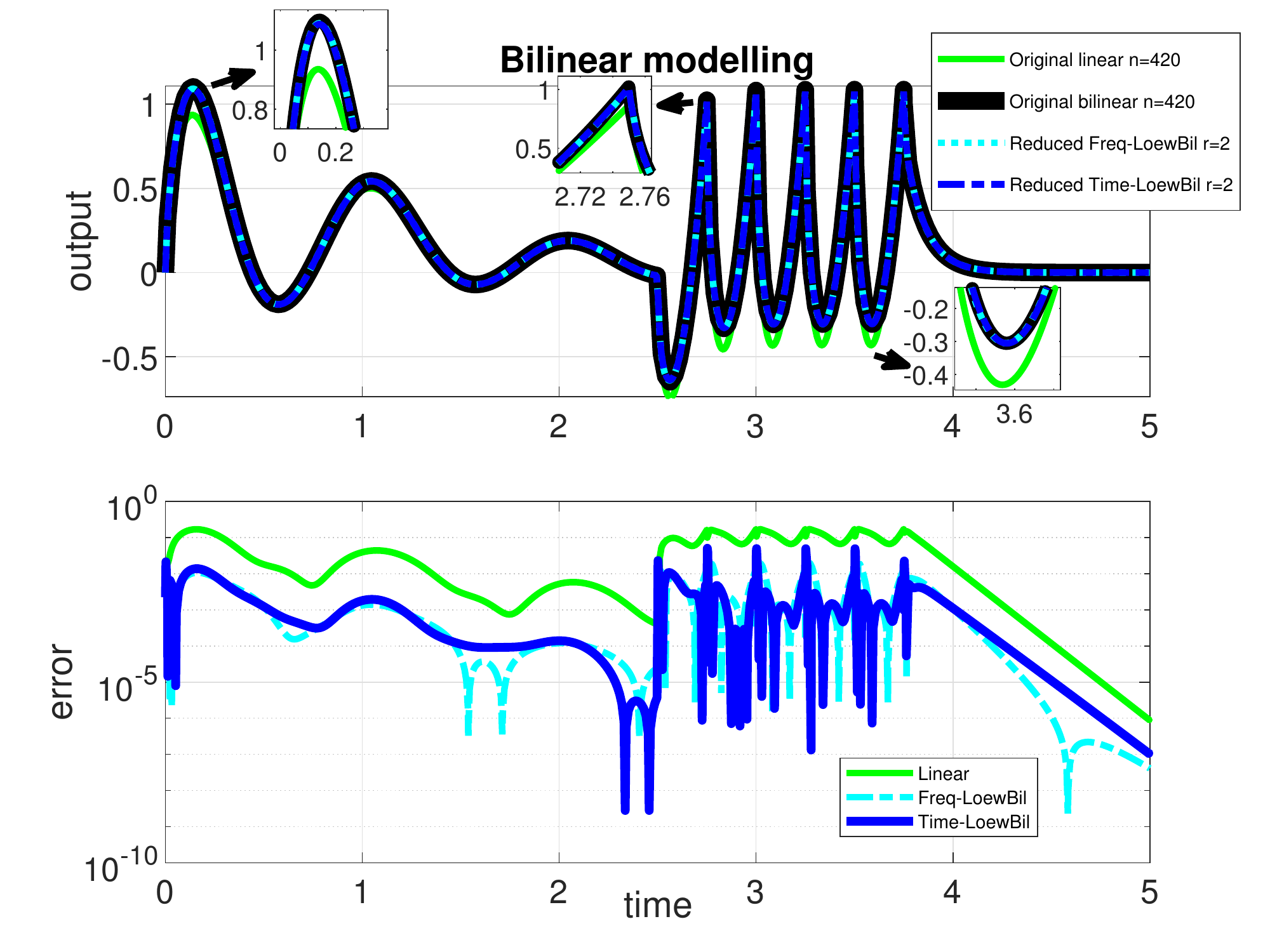}
\caption{Time domain simulation for the Burgers' equation example; viscosity parameter $\nu$ is set as 1 and the dimension of the semi-discretized model is chosen to be $420$. A comparison among the identified/reduced bilinear of order $r=2$ with the linear and with the frequency domain Loewner bilinear is depicted. The input is chosen as: $u(t)=(1+2\cos(2\pi t))e^{-t},t\in[0,2.5],~u(t)=4\text{sawtooth}(8\pi t),t\in[2.5,3.75],~u(t)=0,t\in[3.75,5]$.}
\label{fig:burgersSimulation}
\end{figure}
\end{example}
\begin{table}[!h]
\caption{Summary of the results from the two examples with Time-LoewBil and comparison with \cite{AGIbil} for the Burgers' example 2 of dimension $n=420$.}\label{tab:5}   
\begin{tabular}{p{2.5cm}p{2cm}p{3.5cm}p{3cm}}
\hline\noalign{\smallskip}
Error quantification & Example 1 & Example 2 Time-LoewBil & Example 2 Freq-LoewBil \\
\noalign{\smallskip}\svhline\noalign{\smallskip}
$\Vert H_{1}-\tilde{H_{1}} \Vert_{\text{max}}$ & $5.077\cdot 10^{-3}$  & $2.937\cdot 10^{-3}$ & $4.430\cdot 10^{-3}$ \\[1mm]
$\Vert y(t)-y_{l}(t) \Vert_{\text{max}}$ & $1.213\cdot 10^{-1}$  & $1.699\cdot 10^{-1}$ & $1.699\cdot 10^{-1}$\\[1mm]
$\Vert H_{2}-\tilde{H_{2}} \Vert_{\text{max}}$ & $2.794\cdot 10^{-2}$  & $3.077\cdot 10^{-3}$ & $2.991\cdot 10^{-3}$\\[1mm]
$\Vert y(t)-\tilde{y}_{b}(t) \Vert_{\text{max}}$ & $\bf{2.739\cdot 10^{-4}}$  & $\bf{5.032\cdot 10^{-2}}$ & $\bf{5.278\cdot 10^{-2}}$ \\[1mm]
\noalign{\smallskip}\hline\noalign{\smallskip}
\end{tabular}\\
$^d$ The evaluations of the kernels and the outputs ($y_l$: linear, $\tilde{y}_b$ reduced bilinear (r=2)) took place over the domains depicted in Figs.\;(\ref{fig:kernels},\ref{fig:evalut},\ref{fig:burgersKernels},\ref{fig:burgersSimulation}).
\end{table}

\begin{remark}[Computational cost for the discretized Burgers' model of dimension $420$]\\
	The proposed time-domain Loewner bilinear method uses measurements corresponding to symmetric transfer functions. Such values can be directly inferred from time-domain data by processing the spectral domain, i.e., by computing the FFT of the observed output signals for oscillatory input signals. All experiments
		were performed on a computer with 12 GB RAM and an Intel(R) Core(TM) i7-10510U CPU
		running at 1.80 GHz, 2304 Mhz, 4 Cores, 8 Logical Processors. To simulate a system of dimension $420$, each measurement took $\sim 3$min. {So, the data acquisition cost was reported in the range of 1 or 2 hours where the identification/reduction part was almost direct.} The proposed method seems to efficiently for moderate dimensions; for large-scale problems, the computational issues that appear belong to the class of "embarrassingly parallel" tasks; as the simulations are independent to each other, one can easily speed up the whole process by using instead parallel clusters.
\end{remark}

\begin{remark}[Discussion and comparison between the two methods]
	In what follows we will state the pluses and minuses of the two methods applied for the second numerical example.
\paragraph{{The frequency Loewner bilinear framework (Freq-LoewBil)}}
\begin{itemize}
	\item Pluses: recovers the original bilinear system with high accuracy, incorporates linear and nonlinear transfer
	function measurements in a coupled way ("all at once"), can be easily extended to cope with
	higher-order regular kernels, can also be viewed
	as a Petrov-Galerkin projection-based moment matching approach.
	\item Minuses: It is not completely clear how to measure/obtain the frequency-domain data needed for this method; it uses measurements of regular transfer functions which cannot be (directly) inferred from time-domain simulations.
\end{itemize}
\paragraph{{The time Loewner bilinear framework (Time-LoewBil)}}
\begin{itemize}
	\item Pluses: It uses measurements corresponding to symmetric transfer functions. Such values
	can be directly inferred from time-domain data by processing the spectral domain, i.e., by
	computing the FFT of the observed output signals for oscillatory input signals.
	\item Minuses: The fitted bilinear model is as good as the fitted linear model (it relies on the linear fit). As opposed to the first method, it fits the linear and
	nonlinear parts separately (not "all at once"). It introduces additional errors due to conversion
	from the time domain to the frequency domain. The latter disadvantage could also occur for the
	method in \cite{AGIbil}, provided that "regular transfer function" measurements could be successfully inferred
	from time-domain data.
\end{itemize}
\end{remark}

\vspace{-8mm}

\section{Conclusion}

 The proposed method offers approximate bilinear system identification from time-domain measurements, since it is not possible to measure the corresponding kernels exactly. Our proposed method uses only \textit{input-output} measurements without requiring state-space access. What makes this algorithm feasible is the combination of the data-driven Loewner framework with the nonlinear Volterra series framework.

We have shown that for the noise free case, the proposed method achieves system identification from time-domain measurements through the symmetric kernels. Further study is required to quantify the effects of the noise introduced  by the truncation of the Volterra series (in the $\ell$-stage approximation). {All the time-domain numerical simulations have been implemented by means of the backward Euler approximation scheme which certifies that this method can handle some level of numerical noise. Multi-stepping methods, e.g., Runge Kutta can offer a significance improvement to the results and reduce the influence of numerical noise.}

The variational approach is a theoretical method to identify regular kernels which are appropriate for system identification purposes \cite{Rug81}. However, these kernels do not have a physical meaning, i.e., cannot be directly measured from time-domain simulations. This is not an issue for the growing exponential approach. The derived transfer functions by means of this method can be measured from time-domain data. The difficulty in combining both derivations, i.e., symmetric  and regular is also explained from the $n^{th}$ dimensional integral that connects those through the triangular kernels.
Extensions to the MIMO case and to other nonlinearity structures, e.g., quadratic or bilinear-quadratic etc., are promising endeavors that will be the matter of future research.

\scriptsize{
 \bibliographystyle{spmpsci}
 \bibliography{test}}
\end{document}